\documentclass[11pt]{amsart}
\pdfoutput=1
\usepackage{latexsym, graphicx, epsfig, amsmath, amsfonts,amssymb}
\usepackage{multirow}
 \usepackage{diagbox}
\usepackage{color}
 \usepackage[percent]{overpic}
\usepackage{algorithm}
\usepackage{algorithmicx}
\usepackage{algpseudocode}
\floatname{algorithm}{Algorithm}
\usepackage{threeparttable}
\usepackage{booktabs}

\usepackage{epstopdf,mathtools,bbm}
\usepackage[left=1.25in,right=1.25in,top=1in]{geometry}

\def\ds{\displaystyle}
\def\mb{\mathbf}
\def\R{\mathbb{R}}

\newcounter{Rownumber}

\title{Stein variational gradient descent with local approximations}

\author{Liang Yan}
\thanks{School of Mathematics, Southeast University, Nanjing, China. Email: yanliang@seu.edu.cn. L. Yan is supported by NSF of China (No.11771081), the science challenge project (No. TZ2018001) and  Zhishan Young Scholar Program of SEU}

\author{Tao Zhou}
\thanks{LSEC, Institute of Computational Mathematics, Academy of Mathematics and Systems
Science, Chinese Academy of Sciences, Beijing 100190, China. Email: tzhou@lsec.cc.ac.cn. T. Zhou is partially supported  by the National Key R$\&$D Program of China(No. 2020YFA0712000), the NSF of China (under grant numbers 11822111, 11688101and 11731006), the science challenge project (No. TZ2018001),  the Strategic Priority Research Program of Chinese Academy of Sciences (No.  XDA25000404) and youth innovation promotion association (CAS)}

\date{October 10, 2020}
\begin{document}

\begin{abstract}
Bayesian computation plays  an important role in modern machine learning and statistics to reason about uncertainty. A key computational challenge in Bayesian inference is to develop efficient techniques to approximate, or draw samples from posterior distributions.  Stein variational gradient decent (SVGD) has been shown to be a powerful approximate inference algorithm for this issue.  However, the vanilla SVGD requires calculating the gradient of the target density and cannot be applied when the gradient is unavailable or too expensive to evaluate.  In this paper we explore one way to address this challenge by the construction of a local surrogate for the target distribution in which the gradient can be obtained in a much more computationally feasible manner. More specifically, we approximate the forward model using a deep neural network (DNN)  which is trained on a carefully chosen training set, which also determines the quality of the surrogate.  To this end, we propose a general adaptation procedure to refine the local approximation  online without destroying the convergence of the resulting SVGD. This significantly reduces the computational cost of SVGD and leads to a suite of algorithms that are straightforward to implement.  The new algorithm is illustrated on a set of challenging Bayesian inverse problems, and numerical experiments demonstrate a clear improvement in performance and applicability of standard SVGD.
\end{abstract}

%\begin{keywords}
%Stein variational gradient descent, Bayesian inversion, surrogate modeling, deep neural network.
%\end{keywords}

\pagestyle{myheadings}
\thispagestyle{plain}

\maketitle

%%%% Start %%%%%%
\section{Introduction}

Bayesian inference lies at the heart of many machine learning models in both academia and industry \cite{Bishop2006pattern}. It provides a powerful framework for  uncertainty quantification of various statistical learning models.   The central tasks of Bayesian inference are to compute the {\it posterior distribution} of certain unknown parameters and evaluate some statistics, e.g., expectation, variance etc, of a given quantity of interest \cite{Kaipio+Somersalo2005,Stuart2010}.  The posterior is typically not of analytical form and cannot be easily interrogated.  Thus, it is of major importance to develop efficient approximation techniques to tackle the intractable posterior distribution  that arise in Bayesian inference and prediction problems\cite{Gelman2013bayesian}. 

Over the past decades, large amounts of methods have been proposed to approximate the posterior distribution, including  Markov Chain Monte Carlo (MCMC)  \cite{Brooks2011} and variational inference (VI) \cite{Blei2017variational}.  MCMC methods work by simulating Markov chains whose stationary distributions match the posterior distributions.  Despite  nice asymptotic  theoretical properties, MCMC methods are widely criticized for their slow convergence rate in practice. In complex posterior distribution, the samples from MCMC are often found to have high auto-correlation across time, meaning that the Markov chains explore very slowly in the configuration space.   To address this challenge, several efficient strategies that use  geometry information of the posterior have been developed\cite{Bardsley2014SISC,Lan2016emulation,Martin2012SNMC}.  However,  many exact MCMC methods are computing expensive, and recent years have seen the introduction of surrogate-based MCMC procedures (see, e.g., \cite{Conrad2016JASA,Frangos+Marzouk+Willcox2010,stuart+teckentrup2016,yan+guo2015,Yan+Zhang2017IP,Yan+Zhou19JCP, Yan+Zhou2019ADNN}) to increase sampling speed.  VI methods aim to approximate the posterior by a tractable variational distribution. They transform the inference problem as an optimization problem, which minimizes some kind of distance functional over a prescribed family of known distributions \cite{Blei2017variational}.  Although variational techniques enjoy faster computations, the class of approximations used is often limited, e.g., mean-field approximations, implying that no solution is ever able to resemble the true posterior. This is a widely raised objection to variational methods, in that unlike MCMC, even in the asymptotic regime we are unable recover the true posterior.

Recently, a set of particle-based VI methods \cite{Chen2019projected,Detommaso2018stein,Garbuno2020interacting,Han2018stein,Li2020stochastic,Liu2017riemannian,Liu2017stein,Liu2016stein,Lu2019scaling,Wang2019stein} have been proposed to bridge the gap between VI and MCMC techniques. Those methods use a certain number of samples, or particles, to represent the approximating distribution (like MCMC), and update the particles by  solving an optimization problem (like VI).  They have greater non-parametric flexibility than VI, and are also more particle-efficient than MCMC, since they make full use of particles by taking particle interaction into account.   One such approach, known as Stein Variational Gradient Descent (SVGD)\cite{Liu2016stein}, has attracted much attention.  SVGD is a deterministic sampling algorithm that iteratively transports a set of particles to approximate given distributions, based on a gradient-based update that guarantees to optimally decrease the KL divergence within a function space\cite{Liu2017stein}. SVGD has been shown to provide a fast and flexible alternative to traditional methods such as MCMC and VI in various challenging applications.   Despite these successes and achievements,  there are a number of drawbacks of SVGD that limit their power and hamper their wider adoption as a default method for Bayesian inference. It is one of these limitations, the need of the gradient calculation, that we address in this paper. 

SVGD  requires repeated gradient calculations that involve the whole observed data and maybe a complicated model.
When the gradient information of the target distribution is not  available in closed-form or too expensive to evaluate, standard SVGD methods no longer apply. In fact, the posterior distribution of many Bayesian inference problems only available as a black-box density function and the gradient cannot be calculated analytically.   Recently, Han and Liu \cite{Han2018stein}  develop a gradient-free variant of SVGD(GF-SVGD), which replaces the true gradient with a surrogate gradient, and corrects the induced bias using an importance weight. But this re-weighting approach limits the potential scalability of GF-SVGD  since it requires evaluation of the forward model for each particle per SVGD iteration; when the model is computationally expensive, the GF-SVGD can quickly become computationally prohibitive.  A natural question is how to construct an effective approximation of the gradient information that provides a good balance between accuracy and computation cost.

In this work, we address these challenging problems by  constructing a local surrogate that provides effective approximation of the forward model using deep neural network (DNN). To produce the training sets used for training the  DNN approximations, we will introduce a general adaptation framework  to refine the training set.   In summary, we leverage the advances in different fields of SVGD and emulation, in order to design algorithms which build (a) better surrogate models (a.k.a  {\it emulator}) and (b) more efficient SVGD algorithms. The novel algorithms are adaptive schemes which automatically select the nodes to train the DNN and of the resulting  emulator. Namely, the set of the training points used by the emulator is sequentially  updates by solving a simply optimization problem. 

The organization of the paper  is as follows. In the next section, we review the  Bayesian
inversion and the original  SVGD.  Then we introduce the details of our new method, namely SVGD with local approximations in Section 3. Section 4 includes experiments to validate our method and comparisons with standard SVGD methods on serval numerical examples. Finally, we conclude the paper with Section 5. 

\section {Backgrounds}\label{sec:setup}
In this section, we first give a brief overview of the Bayesian inference in Section \ref{sec:2.1}. Then we will introduce the main setup of the Stein variational gradient descent (SVGD) in Section \ref{sec:2.2}.
\subsection{Bayesian inference}\label{sec:2.1}
The goal of Bayesian inference is to compute the posterior distribution $p(x|y)$ of the unknown parameter $x \in X \subseteq \R^d$ conditioned on the data $y\in Y \subseteq \R^n$. In principle, this posterior may be computed from the likelihood function $\mathcal{L}(x|y, \mb{f}) =p(y|x)$ and the prior $p_0(x)$ using the Bayes' rule
$$
\pi(x):=p(x|y) = \frac{\mathcal{L}(x|y, \mb{f}) p_0(x)}{\int \mathcal{L}(x|y, \mb{f}) p_0(x)dx}.$$
Here we denote the posterior as $\pi(x)$ for convenience of notation,  and $\mb{f}: X \rightarrow Y$ is a parameter-to-observable map.  The forward model $\mb{f}$ may enter the likelihood function in various ways. For instance, if $y=\mb{f}(x)+\eta$, where $\eta \sim p_{\eta}$ represents some measurement of model error, then $\mathcal{L}(x|y, \mb{f}) = p_{\eta}(y-\mb{f}(x))$.

Except for reduced classes of distributions like conjugate priors, the integral in the denominator is usually intractable; no general explicit expressions of the posterior are available. Thus, several techniques have been proposed to perform approximate posterior inference. In this work, we consider one popular particle variational inference technique, Stein variational gradient descent (SVGD) to approximate sampling. 

\subsection{Stein Variational Gradient Descent (SVGD)}\label{sec:2.2}

The goal of SVGD is to find a set of particles $\mb{X}=\{x_i\}^N_{i=1}$  to approximate the target distribution  $\pi$, such that the empirical distribution $q(x) =\frac{1}{N}\sum^N_{i=1}  \delta (x-x_i)$ of the particles weakly converges to $\pi$ when $N$ is large. To achieve this, we initialize the particles with some simple distribution, e.g., the prior $p_0(x)$, and iteratively updating them  with a deterministic transformation of form
\begin{equation}\label{iter_eq}
x_{i} \leftarrow x_{i}+\epsilon \phi\left(x_{i}\right), \quad \forall i=1, \ldots, N,
\end{equation}
where $\epsilon$ is a step size, and $\phi: \R^d \rightarrow \R^d$ is a perturbation direction, or  velocity field, which should be chosen to maximumly decrease the KL divergence between the distribution of particles and the target distribution $\pi$. The optimal choice of $\phi$ can be transformed into the following functional optimization problem \cite{Liu2016stein}:
\begin{equation}\label{opt_eq}
\phi^{*}=\underset{\phi \in \mathcal{S}}{\arg \max }\left\{-\left.\frac{d}{d \epsilon} \mathrm{KL}\left(q_{[\epsilon \phi]} \| \pi\right)\right|_{\epsilon=0}\right\},
\end{equation}
where $\mathcal{S}$ is the set of candidate perturbation direction that we optimize over,  $q_{[\epsilon \phi]}$ denotes the distribution of update particles $x' =x+\epsilon \phi(x)$. 

Critically, the gradient of KL divergence in (\ref{opt_eq}) equals a simple linear functional of $\phi$, allowing us to obtain a closed form solution for the optimal $\phi$.  Liu and Wang \cite{Liu2016stein} showed that 
\begin{equation}\label{stein_eq}
%\begin{array}{r}
-\left. \frac{d}{d \epsilon} \mathrm{KL}\left(q_{[\epsilon \phi]} \| \pi\right)\right|_{\epsilon=0}=\mathbb{E}_{x \sim q}\left[\mathcal{A}_{\pi}^{\top} \phi(x)\right],\, \text { with } \mathcal{A}_{\pi}^{\top} \phi(x)=\nabla_{x} \log\pi(x)^{\top} \phi(x)+\nabla_{x}^{\top} \phi(x),
%\end{array}
\end{equation}
where $\mathcal{A}_{\pi}$ is a differential operator called {\it Stein operator} and is formally viewed as a column vector similar to the gradient operate $\nabla_{x}$.  In SVGD, we often choose $\mathcal{S}$ to be the unit ball of a vector-valued reproducing kernel Hilbert space(RKHS) $\mathcal{H}^d = \mathcal{H} \times \cdots \times \mathcal{H}$ with each $\mathcal{H}$ associating with a positive definite kernel $\kappa(x,x')$, that is, $\mathcal{S} =\{\phi \in \mathcal{H}^d : \|\phi\|_{\mathcal{H}^d} \leq 1\}$.   With (\ref{stein_eq}), it was shown in  \cite{Liu2016stein} that the solution of (\ref{opt_eq}) has a closed form expression given by 

\begin{equation}\label{ker_eq}
\phi^{*}(\cdot) \propto \mathbb{E}_{x\sim q}[\mathcal{A}_{\pi} \kappa(x, \cdot)]=\mathbb{E}_{x \sim q}\left[\nabla_{x} \log \pi(x) \kappa(x, \cdot)+\nabla_{x} \kappa(x, \cdot)\right].
\end{equation}
Such $\phi^{*}$ provides the best update direction for the particles within RKHS $\mathcal{H}^d$.  

There are several ways to approximate the expectation in (\ref{ker_eq}). For instance, by taking $q$ to be the empirical measure of the particles, i.e., $q(x)=\frac{1}{N}\sum^N_{i=1} \delta(x-x_i)$, we obtain
\begin{equation}\label{appker_eq}
\phi^{*}(z) \approx \frac{1}{N}\sum^N_{j=1}\left[\nabla_{x_j} \log \pi(x_j) \kappa(x_j, z)+\nabla_{x_j} \kappa(x_j, z)\right].
\end{equation}
Using Eqs. (\ref{iter_eq}) and (\ref{appker_eq}), we obtain the SVGD algorithm which is summarized in Algorithm \ref{alg:svgd}. The $\epsilon_l$ in line 3 is used as a sequence of  step size. For some applications, one may simply set $\epsilon_l =\epsilon$ to be a constant and gets reasonably good results. However, setting $\epsilon_l$ as a constant may yield divergent sequences in many high dimensional problems\cite{Robbins1951stochastic}. One may decrease $\epsilon_l$ to obtain  convergent sequences. In this work, we use the \texttt{Ada-Grad} \cite{Liu2016stein,Zeiler2012Adadelta} approach to choose $\epsilon_l$.

\begin{algorithm}[t]  
  \caption{Stein variational gradient descent (SVGD)}  
  \label{alg:svgd}  
  \begin{algorithmic}[1]  
   \Procedure{RunSvgd}{$\mb{X}^{(1)}, \mathcal{L},\mb{f}, y, p_0, \kappa,T$}
\For {$l=1,\cdots, T$}
  \State  $x^{(l+1)}_i \leftarrow  x^{(l)}_i+\epsilon_l Q_l(x^{(l)}_i), \,\, \forall i=1,\ldots, N,$  where
\begin{equation}\label{svgdeq}
\begin{aligned}
  &Q_l(x^{(l)}_i)=\frac{1}{N}\sum^N_{j=1}\left[\nabla_{x^{(l)}_j} \log \pi(x^{(l)}_j) \kappa(x^{(l)}_j, x^{(l)}_i)+\nabla_{x^{(l)}_j} \kappa(x^{(l)}_j, x^{(l)}_i)\right], \,\, \text{where} \\
    &\nabla_{x} \log \pi(x) =\nabla_{x} \Big(\log \mathcal{L}(x|y,\mb{f})+\log p_0(x)\Big)
    \end{aligned}
\end{equation}
%\begin{equation}\label{svgdeq}
%  Q_l(x^{(l)}_i)=\frac{1}{N}\sum^N_{j=1}\left[\nabla_{x^{(l)}_j} \log \pi(x^{(l)}_j) \kappa(x^{(l)}_j, x^{(l)}_i)+\nabla_{x^{(l)}_j} \kappa(x^{(l)}_j, x^{(l)}_i)\right], 
%\end{equation}
% and the gradient of the log-target density is defined as 
% $$
% \nabla_{x} \log \pi(x) =\nabla_{x} \Big(\log \mathcal{L}(x|y,\mb{f})+\log p_0(x)\Big).
% $$
 \EndFor 
 \EndProcedure
  \end{algorithmic}  
\end{algorithm}

Notice that the perturbation direction of SVGD, i.e. $Q_l$ in Algorithm \ref{alg:svgd}  (Eq.(\ref{svgdeq})) has two terms. The first term corresponds to a weighted average steepest descent direction of the log-target density. This term is responsible for transporting particles towards high-probability regions of $\pi$. In contrast, the second term can be viewed as a  {\it repulsion force} that spreads the particles along the support of $\pi$, preventing them from collapsing around the mode of $\pi$. Therefore, there are two problems that must be addressed to successfully use the SVGD: 1) efficient computation of the derivatives of the expected log-target density $\pi(x)$, and 2) the efficiency choice of reproducing kernel $\kappa$.  The first problem is focus of this paper.  In next Section, we will introduce a framework with a local approximation to address this issue.

To address the second problem, we focussed on the Gaussian kernel 
\begin{equation}\label{RBFker}
\kappa(x,x') =\exp\Big(-\frac{1}{h}\|x-x'\|^2_2\Big),
\end{equation}
 for some length-scale parameter $h$. This  formula is used in the original analysis of Wang and Liu \cite{Liu2016stein}. In order to improve the performance of the algorithm, several {\it pre-conditioning} of the kernel were recently proposed in context of SVGD, see, e.g.  \cite{Wang2020particle} and references therein.  In this paper, we only focus on the Gaussian kernel with a fixed bandwidth. Notice that,  our framework present in next section can easily extend to the pre-conditioning cases, but it goes beyond the scope of this paper.
\section{SVGD with local approximations}\label{sec:method}
The standard SVGD requires the gradient of the target $\pi$ (or the gradient of the forward model $\mb{f}$). Unfortunately, the gradient information of the target distribution is not always available in practice. In some cases, the forward model of interest is only available as a black-box function and the gradient cannot be calculated analytically; in other cases, it may be computationally too expensive to calculate the gradient, in particular, the systems modeled by partial differential equations (PDEs).  To address these challenges,  we start with introducing our framework for SVGD based on local approximations.

\subsection{SVGD with local approximation}

We assume that the gradient information of forward model is not available; or the forward model and the corresponding  gradient evaluation are computationally expensive --- requiring, for example,  high-resolution numerical solutions of  PDEs. In such a setting, the computational bottleneck of SVGD is dominated by the cost of the gradient evaluations required by $Q_l$.  Motivated  by recently fast MCMC methods using local approximation \cite{Conrad2016JASA,Yan+Zhou19JCP, Yan+Zhou2019ADNN}, one possible way is to construct a local approximation or ‘surrogate’ of the forward model, and then to deal with the posterior distribution induced by this surrogate.  Specifically, assume that one has a collection of model evaluations, $\mathcal{D}_t: =\{\left(x_i, \mb{f}(x_i)\right)\}^{n_t}_{i=1}$, and a method for constructing an approximation $\tilde{\mb{f}_t}$ of $\mb{f}$ based on those points. Using this approximation $\tilde{\mb{f}_t}$, one can obtain an approximated surrogate posterior
\begin{eqnarray*}\label{ppdf_surrogate}
\widetilde{\pi}(x) \propto \mathcal{L}(x| y,\tilde{\mb{f}_t}) p_0(x).
\end{eqnarray*}
If the evaluation of the gradient of the approximation $\tilde{\mb{f}_t}$ is inexpensive, then the gradient of the log-target density,  $\log \widetilde{\pi}$,  can be evaluated  for a large number of samples, without resorting to additional simulations of the forward model $\mb{f}$ and their derivatives.

It should be noticed that the accuracy of the approximation influences the choice of the training points. Given sufficient training data over the whole prior distribution, the surrogate will be able to accurately approximate the froward model and their gradient. However, this may lose the gained computational efficiency \cite{Yan+Zhou19JCP,Yan+Zhou2019ADNN}.  Notice that our concern in Bayesian inversion  is the posterior distribution. Thus, we have to make sure that the approximation  is accurate enough in the posterior density region  but there is no need to ensure its accuracy everywhere. However, estimation of the high-probability density region is nontrivial as the solution of the Bayesian inference is unknown until the data are available.  In order to improve the accuracy of the approximation,  we need to design an algorithm that is allowed to refine the approximation, as needed, by computing new forward model evaluations in the {\it local} high-probability density region and adding them to the growing training set $\mathcal{D}_t$.

\begin{algorithm}[t]  
  \caption{Sketch of approximate SVGD algorithm}  
  \label{alg:svgdLP}  
  \begin{algorithmic}[1]  
   \Procedure{RunLsvgd}{$\mb{X}^{(1)},\mathcal{D}_1, \mathcal{L},  \mb{f},  \tilde{\mb{f}_1}, y,p_0, \kappa, T, I_{max}$}
\For {$t=1,\cdots, I_{max}$}
 % \State  Compute $\nabla_{x} \log \pi = \nabla_{x} \log (\mathcal{L}(x|\tilde{f},y)p_0(x))$
  \State $\mb{X}^{(t+1)}\leftarrow \Call{RunSvgd}{\mb{X}^{(t)}, \mathcal{L}, \tilde{\mb{f}_t}, y, p_0, \kappa,T}$  
  \State $(\tilde{\mb{f}}_{t+1}, \mathcal{D}_{t+1})\leftarrow \Call{RefineApprox}{\mb{X}^{(t+1)},\mathcal{D}_t,\tilde{\mb{f}_t},\mb{f}}$ 
 \EndFor 
 \EndProcedure
\vspace{0.3cm}

 \Procedure{RefineApprox}{$\mb{X},\mathcal{D},\tilde{\mb{f}},\mb{f}$}
  \State Choose a design point $x^*$ and some indicator criterion using $\mb{X},\tilde{\mb{f}}$ and $\mb{f}$
  \If {approximate needs refinement near $x^*$}
  \State Select new $Q$ sample points $\{x_i\}$ from $\mb{X}$ and grow $\mathcal{D}\leftarrow \mathcal{D} \cup \{x_i, \mb{f}(x_i)\}$. 
  \State Update the approximate $\tilde{\mb{f}}$  using the training set $\mathcal{D}$.
  \EndIf
    \State
  \Return $\tilde{\mb{f}}$  and $\mathcal{D}$

   \EndProcedure
 
  \end{algorithmic}  
\end{algorithm}  

Our approach, outlined in Algorithm 2, is in the same spirit as those previous efforts. Indeed, the sketch in Algorithm 2 is sufficiently general to encompass all the previous efforts mentioned above. The methodology is modular and broken into following steps, each of which can be tackled by different methodologies:
\begin{itemize}
\item  Initialization:  build a  surrogate $\tilde{\mb{f}}_1$ using the initial design pool $\mathcal{D}_1$ with a small size. 
\item Online computations: using the surrogate  $\tilde{\mb{f}}_t$, run the SVGD algorithm $\textsc{RunSvgd}$  for a certain number of steps (say $T$ steps) to obtain the particles $\mb{X}^{(t+1)}$. 
\item Refinement:  choose a design point $x^*$ and an indicator criterion.  Then decide whether the approximation need to be refine nearby of $x^*$. If the surrogate needs refinement, then select new points from the particles  $\mb{X}^{(t+1)}$  to refine the sample set $\mathcal{D}_t$ and  the approximation $\tilde{\mb{f}}_t$. 
\item  Repeated the above procedure for many times (say at most $I_{max}$ times).
\end{itemize}

The remainder of this section expands this outline into a usable algorithm, detailing how to construct the local approximations, when to perform refinement, and how to select new training points to refine the approximations. Intuitively, one can argue that this algorithm will produce accurate samples if $\tilde{\mb{f}}_t$ is close to $\mb{f}$, and that the algorithm will be efficient if the size of $\mathcal{D}_t$ is small and $\tilde{\mb{f}}_t$ is cheap to construct.

\subsection{Deep neural network (DNN) approximation}\label{sec:DNN}
In this section, we describe how to construct a surrogate model.  Here, we only consider  the deep neural networks (DNN) method, and the framework can be extended to other methods such as polynomial chaos expansions\cite{yan+guo2015,Yan+Zhou19JCP} and Gaussian process\cite{stuart+teckentrup2016}.  

The basic idea of using deep neural networks (DNN) to construct a surrogate model is that one can  approximate  an input-output map $\mb{f}$ through a hierarchical abstract layers of latent variables \cite{Goodfellow2016DL}. A typical example is the feedforward neural network, which is also called multi-layer perception (MLP). It consists of a collection of layers that include an input layer, an output layer, and a number of hidden layers.  Specifically, in the $k$th hidden layer, $d_k$ number of neurons are present. Each hidden layer of the network receives an output $z^{(k-1)}\in \R^{d_{k-1}}$ from the previous layer where an affine transformation of the form
\begin{equation}
\mathcal{F}_k(z^{(k-1)})  = \mb{W}^{(k)}z^{(k-1)}+\mb{b}^{(k)},
\end{equation}
is performed. Here  $\mb{W}^{(k)} \in \R^{d_{k}\times d_{k-1}}, \, \mb{b}^{(k)} \in \R^{d_{k}}$ are the weights and biases of the $k$th layer. A nonlinear function $\sigma$ called {\it activation function} is applied to each component of the transformed vector before sending it as an input to the next layer.  Some popular choices for the activation function include \texttt{sigmoid,  hyperbolic tangent,  rectified linear unit (ReLU)},  to name a few \cite{Goodfellow2016DL,Ramachandran2017}. In the current work, we shall use \texttt{Swish} as the activation function \cite{Ramachandran2017,Tripathy+Bilionis2018JCP}.  Using the above notation,  the surrogate using the  neural network representation is given by the composition
\begin{equation}
\tilde{\mb{f}}(x) = (\mathcal{F}_L\circ\sigma\circ\mathcal{F}_{L-1}\circ \cdots \circ\sigma\circ\mathcal{F}_1)(x),
\end{equation}
where the operator $\circ$ is the composition operator, $x=z^{(0)}$ is the input.  

Once the network architecture is defined, one can use optimization tools to find the unknown parameters $\theta =\{ \mb{W}^{(k)},  \mb{b}^{(k)}\}$ based on the training data. Precisely, use the training data set  $\mathcal{D}: =\{(x_i, y_i)\}^{n_t}_{i=1}$, we can define the following minimization problem:
\begin{equation}\label{thetastar}
\arg\min_{\theta} \frac{1}{n_t}\sum^{n_t}_{i=1} \|y_i-\tilde{\mb{f}}(x_i;\theta)\|^2+\beta \Omega(\theta),
\end{equation}
where $\mathcal{J}(\theta; \mathcal{D}) = \frac{1}{n_t}\sum^{n_t}_{i=1} \|y_i-\tilde{\mb{f}}(x_i;\theta)\|^2+\beta \Omega(\theta)$ is the so called loss function, $\Omega(\theta)$ is a regularizer and $\beta$ is  the regularization constant.  For our case, the regularizer $\Omega(\theta) = \|\theta\|^2$.  Solving this problem is generally achieved by the stochastic gradient descent (SGD) algorithm \cite{Bottou2010}.   SGD simply minimizes the function by taking a negative step along an estimate  of the gradient  $\nabla_{\theta} \mathcal{J}(\theta; \tilde{\mathcal{D}} )$ at iteration $k$, where $\tilde{\mathcal{D}} \subseteq \mathcal{D}$ is a small randomly sampled subset of $\mathcal{D}$. The gradients are usually computed through backpropagation.  At each iteration, SGD updates the solution by
\begin{equation}\label{sgditer}
\theta_{k+1}=\theta_k-\lambda \nabla_{\theta} \mathcal{J}(\theta; \tilde{\mathcal{D}} ),
\end{equation}
where $\lambda$ is the learning rate. Recent algorithms that offer adaptive learning rates are available, such as \texttt{Ada-Grad} \cite{Zeiler2012Adadelta}, \texttt{RMSProp} \cite{Tieleman+Hinton2012lecture} and \texttt{Adam} \cite{Kingma2014Adam}, ect.  The present work adopts \texttt{Adam} optimization algorithm.

Notice that the DNN can potentially handle functions with limited regularity and are powerful tools for approximating high dimensional problems. Using DNN to build surrogate models have been gaining great popularity among diverse scientific disciplines \cite{Han+Jentzen+E2018PNAS,Raissi2019JCP,Schwab+Zech2019AA,Tripathy+Bilionis2018JCP,Zhu+Zabaras2018bayesian}.  
It is clear that after obtaining the parameters $\theta$, we have an explicit functional form $\tilde{\mb{f}}(x;\theta)$ and can compute its gradient $\nabla_x\tilde{\mb{f}}(x;\theta)$ easily via the back propagation \cite{Bishop2006pattern}. These approximations can be then substituted into the computation procedure of the SVGD framework.  However, $\mathcal{D}$ is difficult to design in advance. A naive way to choose the training sets is to generate enough data over the whole prior distribution; but a prior based surrogate might not be enough for online computations,  see e.g. \cite{Yan+Zhou19JCP,Yan+Zhou2019RTO}. Thus, in next section, we will present  an adaptive procedure to update the training sets and the DNN model.

\subsection{Refining the local DNN model}
\begin{algorithm}[t]  
  \caption{Refine a local DNN}  
  \label{alg:ldnn}  
  \begin{algorithmic}[1]  
   \Procedure{RefineDnn}{$ \mb{X}^{(t+1)}, \mathcal{D}, \tilde{\mb{f}},\mb{f}, tol, Q, R, \rho$}
 \State Compute $x^{*}$ and the relative error $err(x^*)$ using Eqs. (\ref{mpts}) and  (\ref{conerr}), respectively. 
 \If {$err(x^*)>tol$}
\State $\ell \leftarrow1$
 \While {$\ell \leq Q$}
      \State $x_{\ell } \leftarrow \arg\min_{x' \in \mb{X}^{(t+1)}} \|x'-x^*\|_2,
   \text{subject to} \, \|x' - x\|_2 \geq R,\,\, \forall x \in \mathcal{D}.$
      \If {$x_{\ell}$ is in-existence} {break}
      \Else
      \State $\mathcal{D} \leftarrow \mathcal{D} \cup \{x_{\ell }, \mb{f}(x_{\ell })\}$. 
      \State $\ell  \leftarrow \ell +1$
      \EndIf
   \EndWhile   
%  \State   Select $R$ so that $|\mathcal{B}(x^*,R)| \leq Q$, where
%  $$
%  B(x^*,R):= \{x_i\in \mb{X}^{(t+1)}: \|x_i-x^*\|\leq R, \text{and} \|x_i-x\|\leq R_0, \forall x \in \mathcal{D}\}
%  $$
%  \State $\mathcal{D} \leftarrow \mathcal{D} \cup \{x_i, \mb{f}(x_i)\}$. 
    \EndIf
    \If {$\ell \leq 1$}  
    \State {$R \leftarrow \rho R$} 
    \Else 
     \State  Refine the surrogate model  $\tilde{\mb{f}}$ using $\mathcal{D}$.
       \EndIf
    \State
    \Return $\tilde{\mb{f}}, \mathcal{D}$ and $R$
 \EndProcedure
  \end{algorithmic}  
\end{algorithm}  

As  mentioned previously, the choice of the training points plays a crucial role in training the DNN. First, the quality of the training points directly determines the precision of the associated surrogate. Second, the size of the training sets $\mathcal{D}$ should be controlled for computational efficiency. An ideal set of training points should evenly spread over the posterior density region so that the DNN conditioning on it models the target probability distribution well, however, it is often challenging to do so at the beginning, especially for computation intensive models.  Therefore, it is natural to think of adapting the initial training sets and the surrogate while generating the SVGD.  Starting from some initial training set with a small size, the adaptive method can grow them from some candidate set to a desirable size according to some information criterion.  Considering the setting of SVGD, it is natural to grow the training sets by selecting candidates from previous particles which located in the high probability density regions(HPDR).  On one hand, the SVGD feeds the surrogate with useful candidate points to refine with; on the other hand, the DNN surrogate returns gradient  information efficiently for the SVGD to further explore the parameter space. They form a mutual learning system so that they can learn from each other and gradually improve each other. 

To this end, we first sampling approximate posterior distribution based on the surrogate model $\tilde{\mb{f}}$ for a certain number of iteration steps, namely $T$ steps,  using a standard SVGD algorithm.  The theory and numerical results shown in \cite{Liu2017stein,Liu2016stein} demonstrate that the particles $\mb{X}^{(t+1)}=\{x_i^{(t+1)}\}^N_{i=1}$ obtained by SVGD are in HPDR.  Using those particles, we can compute a posterior  mean point $x^*$ as
\begin{equation}\label{mpts}
x^* = \frac{1}{N} \sum_{i=1}^M x^{(t+1)}_i,
\end{equation}
which is expected to be much closer to the posterior region.  Once we obtain the designed point $x^*$,  we can define an error indicator as 
\begin{equation}\label{conerr}
err(x^*) = \frac{\|\mb{f}(x^*)-\tilde{\mb{f}}(x^*)\|_{2}}{\|\mb{f}(x^*)\|_{2}}.
\end{equation}
If this error indicator exceeds a user-given threshold $tol$, the local DNN need to refine near $x^*$.  Otherwise the error indicator is smaller than $tol$, it means that the surrogate model is still acceptable and we just go ahead.  If refinement of the DNN model at a point $x^*$ is required, we perform the refinement by selecting $Q$ new points $\{x_i\}_{i=1}^Q$ from $\mb{X}^{(t+1)}$,   compute $\mb{f}(x_i)$ and insert the new pair into  $\mathcal{D}$.  In particular, we shall choose $Q$ new points  $\{x_i\}\in \mb{X}^{(t+1)}$ in a ball centered at $x^*$, i.e., 
$x_i \in\mathcal{B}(x^*,R_t):= \{x_i \in\mb{X}^{(t+1)}:\|x_i-x^*\|\leq R_t \}.$  Here the radius $R_t$ is selected to include a fixed number of points $Q$. 

Notice that the strategy of simply adding $x_i$ from $\mb{X}^{(t+1)}$ to $\mathcal{D}$ is inadvisable because it often introduces tightly clustered points.  Instead, a more straightforward and widely used type of experimental design is to choose points in a space-filling fashion. Specifically, we select a new point $x_{i}$ by finding the local minimizer of the problem:
 \begin{equation}\label{Lminer}
 \begin{aligned}
 x_{i}= &\arg\min_{x' \in \mb{X}^{(t+1)}} \|x'-x^*\|_2,\\
 & \text{subject to} \, \|x' - x\|_2 \geq R,\,\, \forall x \in \mathcal{D}.
 \end{aligned}
 \end{equation}
Here $R$ is a constant to control the new added points $\{x_i\}$ are  well separated from the entire set $\mathcal{D}$, and the minimization operator finds a point  lying in the local point set $\mathcal{B}(x^*,R_t)$.  Unfortunately,  the local minimizer $x_{i}$ in (\ref{Lminer}) maybe in-existence.  When it does, we  break the algorithm, and this makes the size of the set $\mathcal{B}(x^*,R_t)$  smaller than $Q$. Specially, while the set $\mathcal{B}(x^*,R_t)$ is empty, meaning that the value of $R$ in Eq. (\ref{Lminer}) is too large and we should use a relatively small $R$. This  motivates us to introduce a constant $0<\rho<1$ to control the value of  $R$. This strategy is detailed  in Algorithm \ref{alg:ldnn}.   After obtaining the new training set  $\mathcal{D}$, we can refine the DNN model $\tilde{\mb{f}}$ in line 17. More specifically,  we initialize the parameter $\theta$ of the DNN  from the pre-trained  model $\tilde{\mb{f}}$ in the online training process. This initialization can be considered as an instance of transfer learning \cite{Yosinski2014transferable}, and we expect a considerable speed-up when solving Eq. (\ref{thetastar}).

\subsection{Algorithm summary}

\begin{algorithm}
  \caption{SVGD with local DNN approximation}  
  \label{alg:LDNN}  
  \begin{algorithmic}[1]  
   \Procedure{RunLsvgd}{$\mb{X}^{(1)},\mathcal{D}_1, \mathcal{L},  \mb{f},  \tilde{\mb{f}_1}, y, p_0, \kappa, T, I_{max}, tol, Q, R, \rho$}
\For {$t=1,\cdots, I_{max}$}
 % \State  Compute $\nabla_{x} \log \pi = \nabla_{x} \log (\mathcal{L}(x|\tilde{f},y)p_0(x))$
  \State $\mb{X}^{(t+1)}\leftarrow \Call{RunSvgd}{\mb{X}^{(t)}, \mathcal{L}, \tilde{\mb{f}_t}, y, p_0, \kappa,T}$  
  \State $(\tilde{\mb{f}}_{t+1}, \mathcal{D}_{t+1}, R)\leftarrow \Call{RefineDNN}{\mb{X}^{(t+1)},\mathcal{D}_t,\tilde{\mb{f}_t},\mb{f}, tol, Q, R, \rho}$ 
 \EndFor 
 \EndProcedure
\vspace{0.1cm}

 \Procedure{RefineDNN}{$\mb{X},\mathcal{D},\tilde{\mb{f}},\mb{f}, tol, Q, R, \rho$}
  \State Compute $x^{*}$ and  $err(x^*)$ using Eqs. (\ref{mpts}) and  (\ref{conerr}), respectively. 
 \If {$err(x^*)>tol$}
\State $\ell \leftarrow 1$
 \While {$\ell \leq Q$}
      \State $x_{\ell } \leftarrow \ds \arg\min_{x' \in \mb{X}^{(t+1)}} \|x'-x^*\|_2,
   \text{subject to} \, \|x' - x\|_2 \geq R,\,\forall x \in \mathcal{D}.$
      \If {$x_{\ell}$ is in-existence} {break}
      \Else
      \State $\mathcal{D} \leftarrow \mathcal{D} \cup \{x_{\ell }, \mb{f}(x_{\ell })\}$. 
      \State $\ell  \leftarrow \ell +1$
      \EndIf
   \EndWhile   
    \EndIf
    \If {$\ell \leq 1$}  
    \State {$R \leftarrow \rho R$} 
    \Else 
     \State  Refine the surrogate model  $\tilde{\mb{f}}$ using $\mathcal{D}$.
       \EndIf
    \State
    \Return $\tilde{\mb{f}}, \mathcal{D}$ and $R$
 \EndProcedure 
  \end{algorithmic}  
\end{algorithm}  

\begin{figure}[h]
\begin{center}
         \begin{overpic}[height=4.4cm,width=3.3cm,trim= 35 10 45 15, clip=true,tics=10]{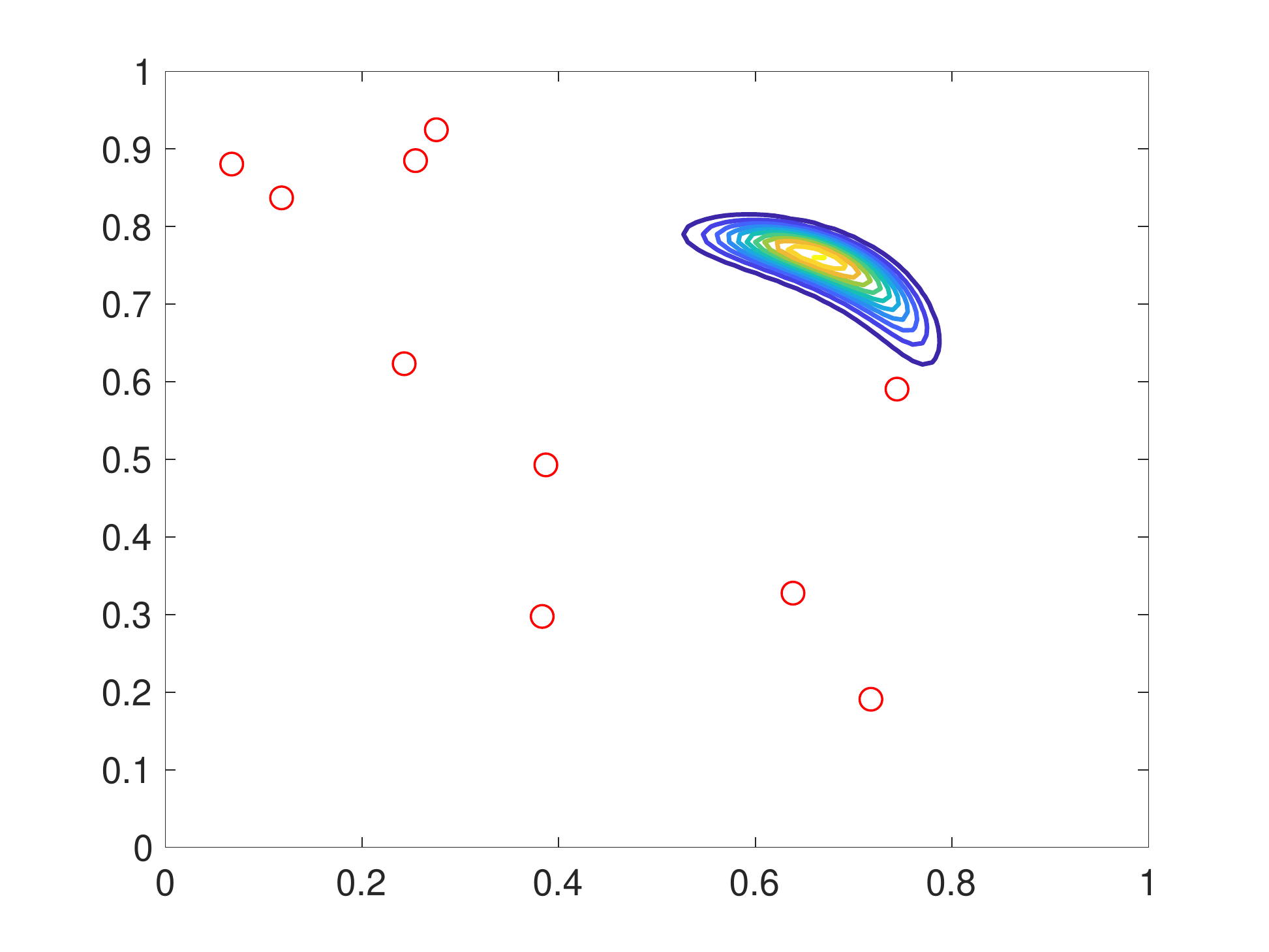}
         \put (5,98) {\scriptsize {\bf Initial design  points}}
          \put (40,12) {\footnotesize \bf Case 1}
        \end{overpic}
    \begin{overpic}[height=4.4cm,width=3.3cm,trim= 35 10 45 15, clip=true,tics=10]{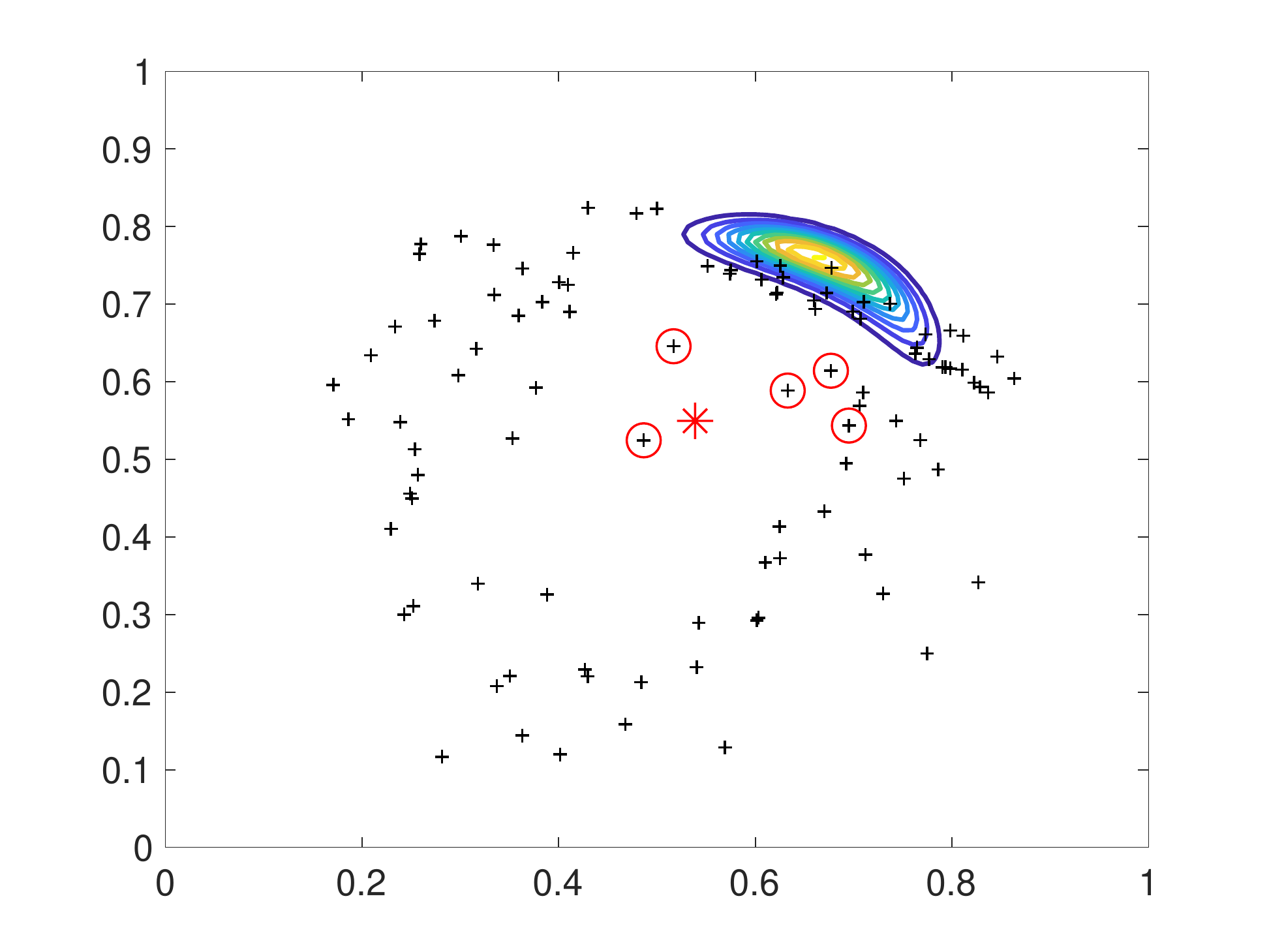}
         \put (18,98) {\scriptsize {\bf First update}}
  \end{overpic}
  \begin{overpic}[height=4.4cm,width=3.3cm,trim=35 10 45 15, clip=true,tics=10]{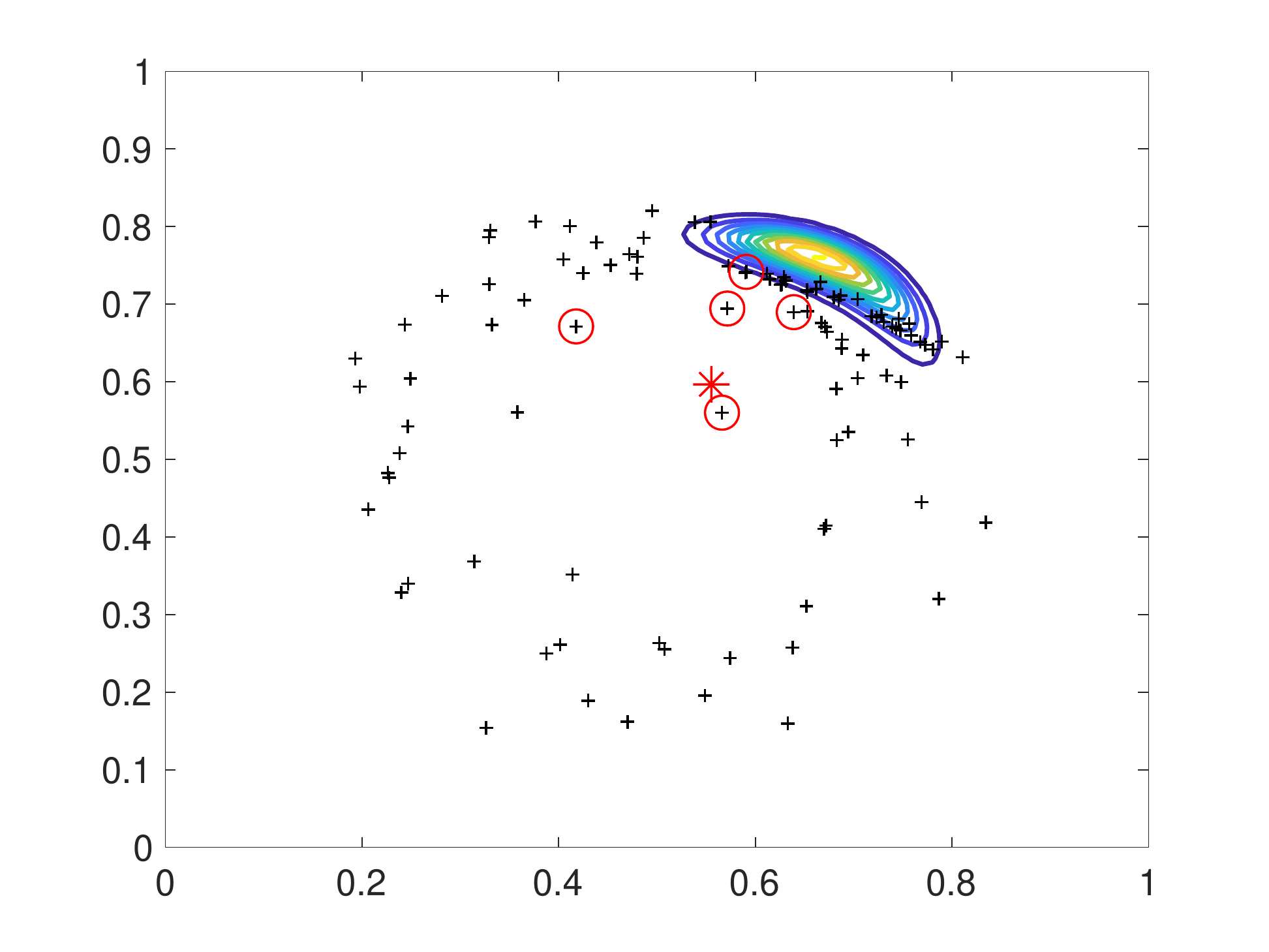}
       \put (18,98) {\scriptsize {\bf Adapting...}}
  \end{overpic}
     \begin{overpic}[height=4.4cm,width=3.3cm,trim= 35 10 45 15, clip=true,tics=10]{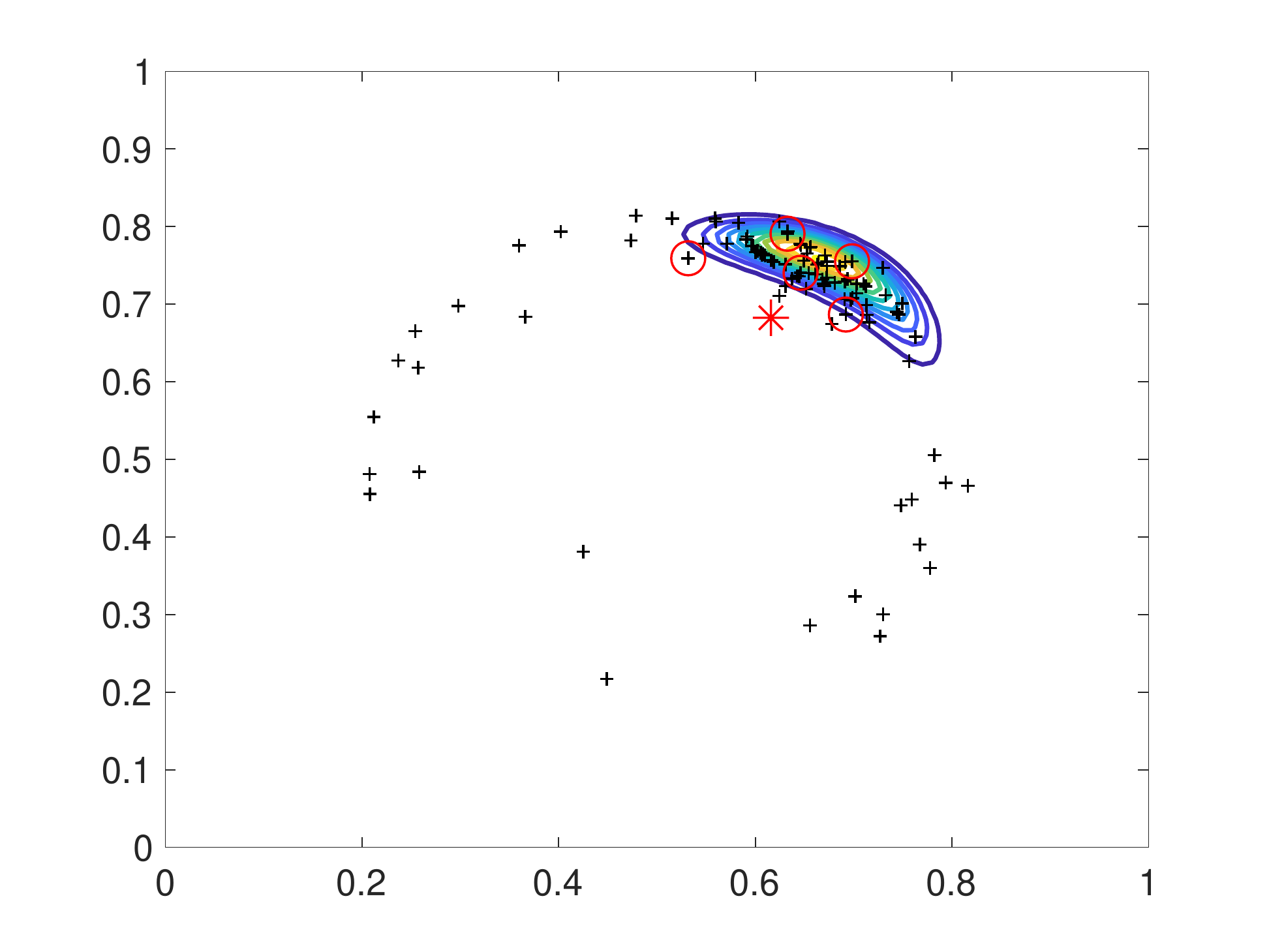}
         \put (18,98) {\scriptsize {\bf Final update}}
  \end{overpic}
     \begin{overpic}[height=4.4cm,width=3.3cm,trim= 35 10 45 15, clip=true,tics=10]{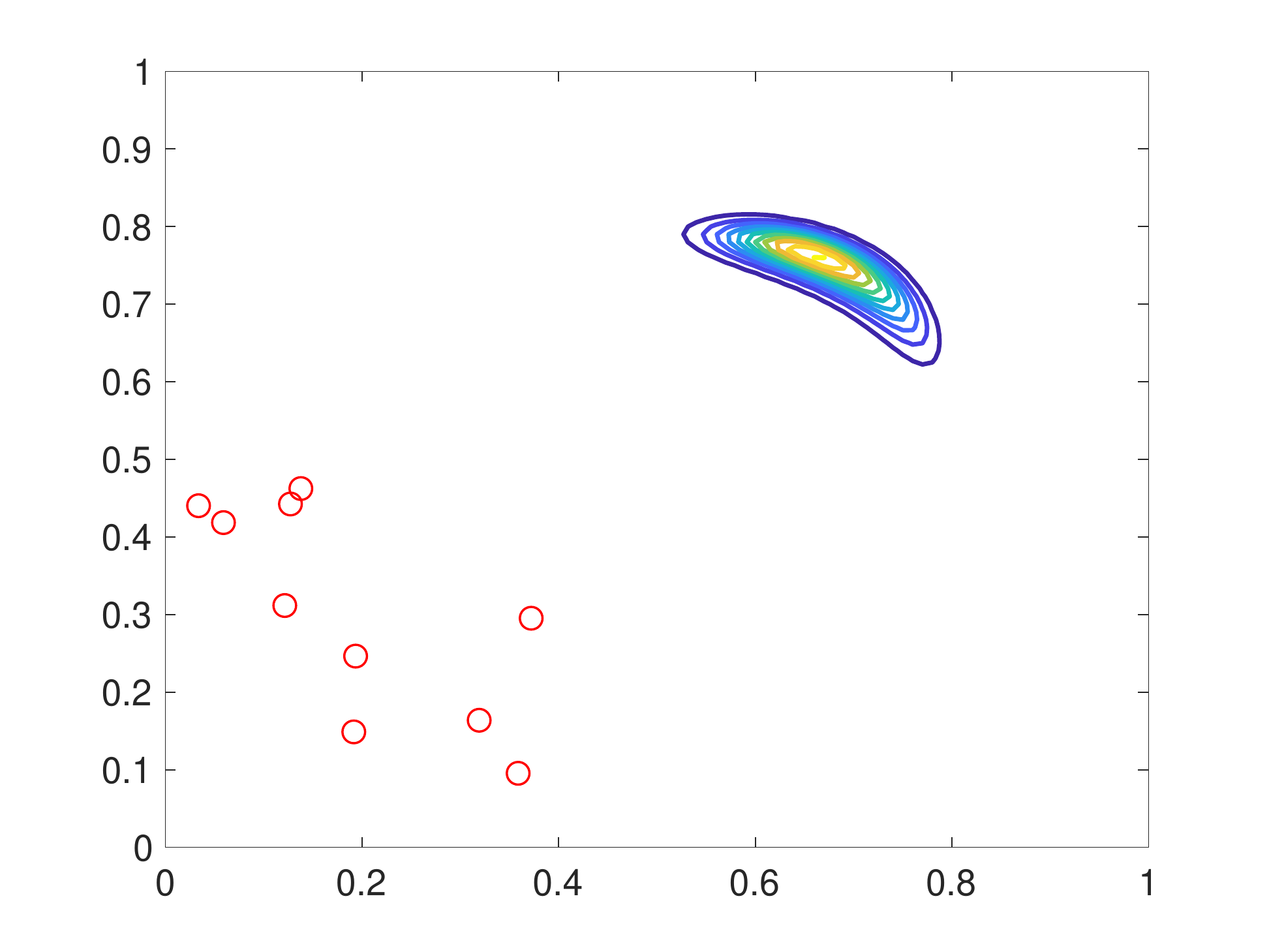}
     \put (40,12) {\footnotesize \bf Case 2}
  \end{overpic}
    \begin{overpic}[height=4.4cm,width=3.3cm,trim= 35 10 45 15, clip=true,tics=10]{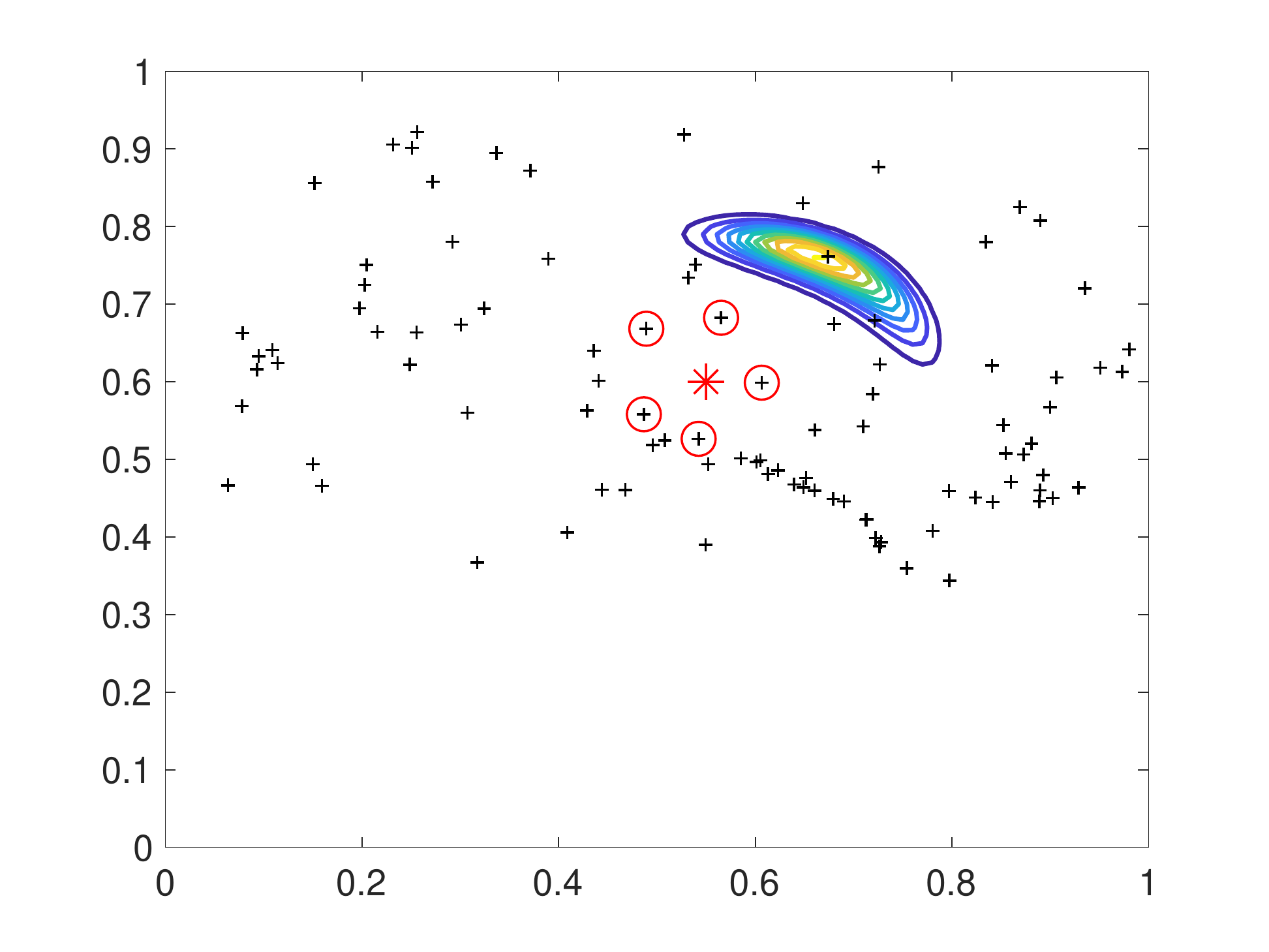}
  \end{overpic}
  \begin{overpic}[height=4.4cm,width=3.3cm,trim=35 10 45 15, clip=true,tics=10]{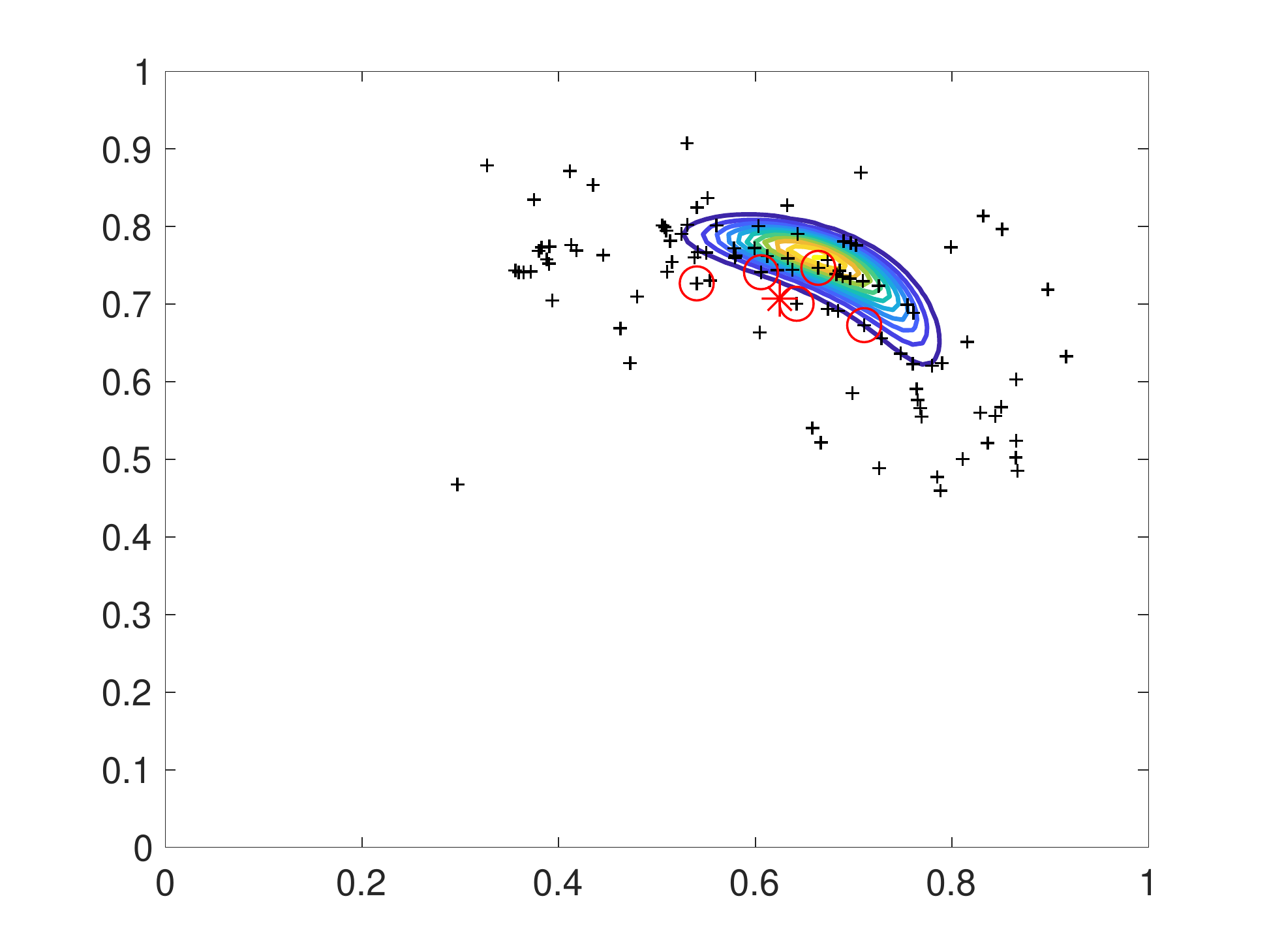}
  \end{overpic}
     \begin{overpic}[height=4.4cm,width=3.3cm,trim= 35 10 45 15, clip=true,tics=10]{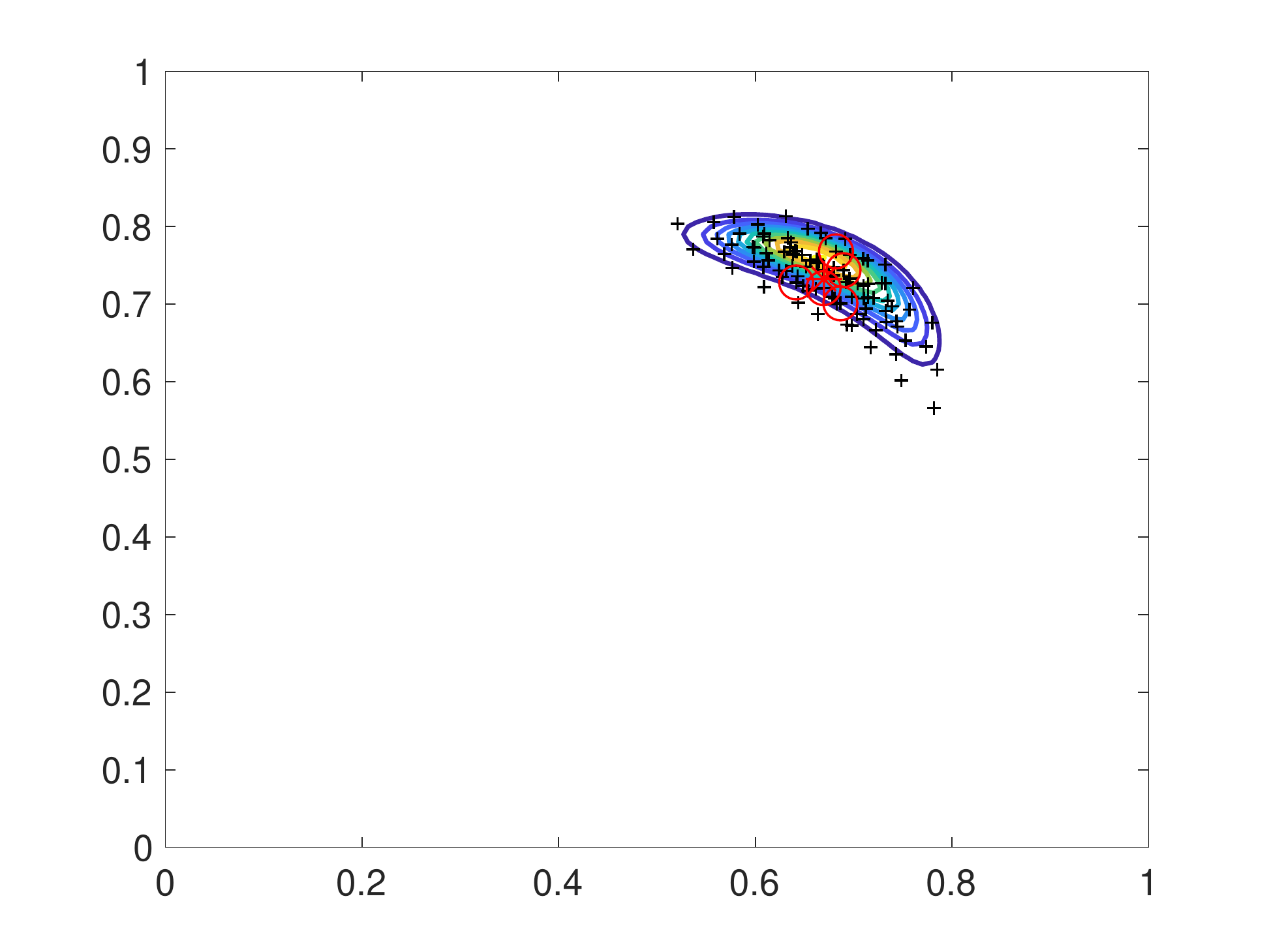} 
  \end{overpic}
  \end{center}
\caption{ The evolution of the training points by the adaptive algorithm. Black $+$ are particles obtained by SVGD, red circles are new design points added to the training set $\mathcal{D}$ and red $*$ is  $x^*$  in each adaptation.  }\label{uppoints_eg2}
  \end{figure}

Our SVGD approach using local DNN approximations (LDNN-based SVGD) is summarized in Algorithm 4. The algorithm proceeds in much the same way as the sketch provided in Algorithm 2.  To summarize, our approach starts from an offline pre-trained DNN surrogate, and then we correct this DNN adaptively using locally training data from the previous particles.   Fig. \ref{uppoints_eg2} shows how the online training points are evolved by the adaptive algorithm. Even starting with a bad design (see the bottom figures), the adaptive algorithm can gradually spread them over the target density region until it reaches the final design based on which DNN accurately emulates the true probability distribution.  More details of this example are shown in Section \ref{heatsource}. 

We review the computational efficiency of Algorithm 4.  Let $q_t$ denote the points in the local ball $\mathcal{B}(x^*,R_t)$ within distance $R_t$. Then the total number of online high-fidelity model evaluations for the LDNN-based SVGD algorithm is $N_{eval} = \sum^{I_{max}}_{t=1} q_t +I_{max} \leq (Q+1) I_{max}$. Note that  this value is independent of the number of particles $N$.  As our empirical results suggest, with appropriate training data good approximation  can be achieved using a moderate  number of training points. On the other hand,  if we run $N_{iter}$ iterations original SVGD algorithm with a  large number of particles $N$ (e.g., $10^2-10^3$), we need  $N_{iter} \times N$ evaluations of the gradient of the high-fidelity model.  Focusing on the number of forward model or the corresponding gradient evaluations, the computational cost of LDNN-based SVGD is  significantly reduced.  Therefore, our proposed method has the potential to provide effective and scalable surrogate-based SVGD that balances accuracy and efficiency well. In addition, the new scheme does not require additional evaluations of the gradient of the exact high-fidelity model, which expands the application of the  original SVGD method. 
\section{Numerical Experiments}\label{sec:tests}

To illustrate the utility of the LDNN-based SVGD, we first apply it to one toy model with analytically tractable posteriors. We then demonstrate the benefits of the local approximation framework for two challenging Bayesian inference problems.  We use those examples to demonstrate the utility of SVGD with local approximation in terms of speed, accuracy of recovery, and probabilistic calibration.  To evaluate how well the particles approximate the posterior, we consider Maximum Mean Discrepancy (MMD)\cite{Gretton2012kernel} related to sample quality. 

To limit the scope, we present a comparison of LDNN-based SVGD to the original SVGD method, as well as to SVGD with global approximation (prior-based DNN).  We refer to these algorithms as LDNN, Direct and DNN, respectively.  In all our numerical tests, the DNN approach was performed with a self-written program which was coded by MATLAB. The optimization procedure is carried out by the \texttt{Adam}  algorithm as mentioned before. The regularization constant in Eq. (\ref{thetastar}) is set to be $\beta=1\times 10^{-6}$, the learning rate is set to be $\lambda= 5\times 10^{-4}$, and the hyper-parameter values of  \texttt{Adam} are chosen based on default recommendations as suggested in \cite{Kingma2014Adam}.   For SVGD, we use the kernel (\ref{RBFker}) with the bandwidth $h= med^2/\log N$, where $med$ is the median of the current $N$ particles. The initial point set for SVGD was chosen by sampling each point independently from the prior $p_0(x)$. The step-size $\epsilon$ for SVGD was set using  \texttt{Ada-Grad},  as in \cite{Liu2016stein}, with master step size 1 and momentum 0.9.  When MMD is applied to evaluate the sample quality, RBF kernel is used and the bandwidth is chosen based on the median distance of the `exact' samples so that all methods use the same bandwidth for a fair comparison.  In the examples below, unless otherwise specified, we use the following parameters  $Q= 5, R=0.2,  tol = 1\times 10^{-2}, \rho =0.8$.   We run the LDNN-based SVGD algorithm for $I_{max} =30$ outer loop iterations, with the number of the inner loop $T=10$.  To make a fair comparison, we run the prior-based DNN and the original  SVGD  for $N_{iter}=I_{max} \times T$ iterations.  All the computations were performed using MATLAB 2018a on an Intel-i7 desktop computer.

\subsection{2D toy example: Double Banana}
In the  first example we consider the ``double-banana" target density initially presented in \cite{Detommaso2018stein}, whose probability density function is 
\begin{equation*}
\pi(x) \propto \exp\Big(-\frac{\|x\|^2_2}{2\sigma_1^2}-\frac{(y-\mb{f}(x))^2}{2\sigma_2^2}\Big).
\end{equation*}
Here $x=[x_1,x_2] \in \R^2$,  $\mb{f}(x)=\log\big((1-x_1)^2+100(x_2-x_1^2)^2\big)$ and $y=\log(30), \sigma_1=1,\sigma_2 =0.3$.  
Though the gradient information of this example can be explicitly calculated, and performing SVGD directly is  at quite low cost, we may still consider whether local approximations can reduce the number of times the exact model must be evaluated. 

\begin{figure}
\begin{center}
         \begin{overpic}[width=0.85\textwidth,trim= 35 10 45 15, clip=true,tics=10]{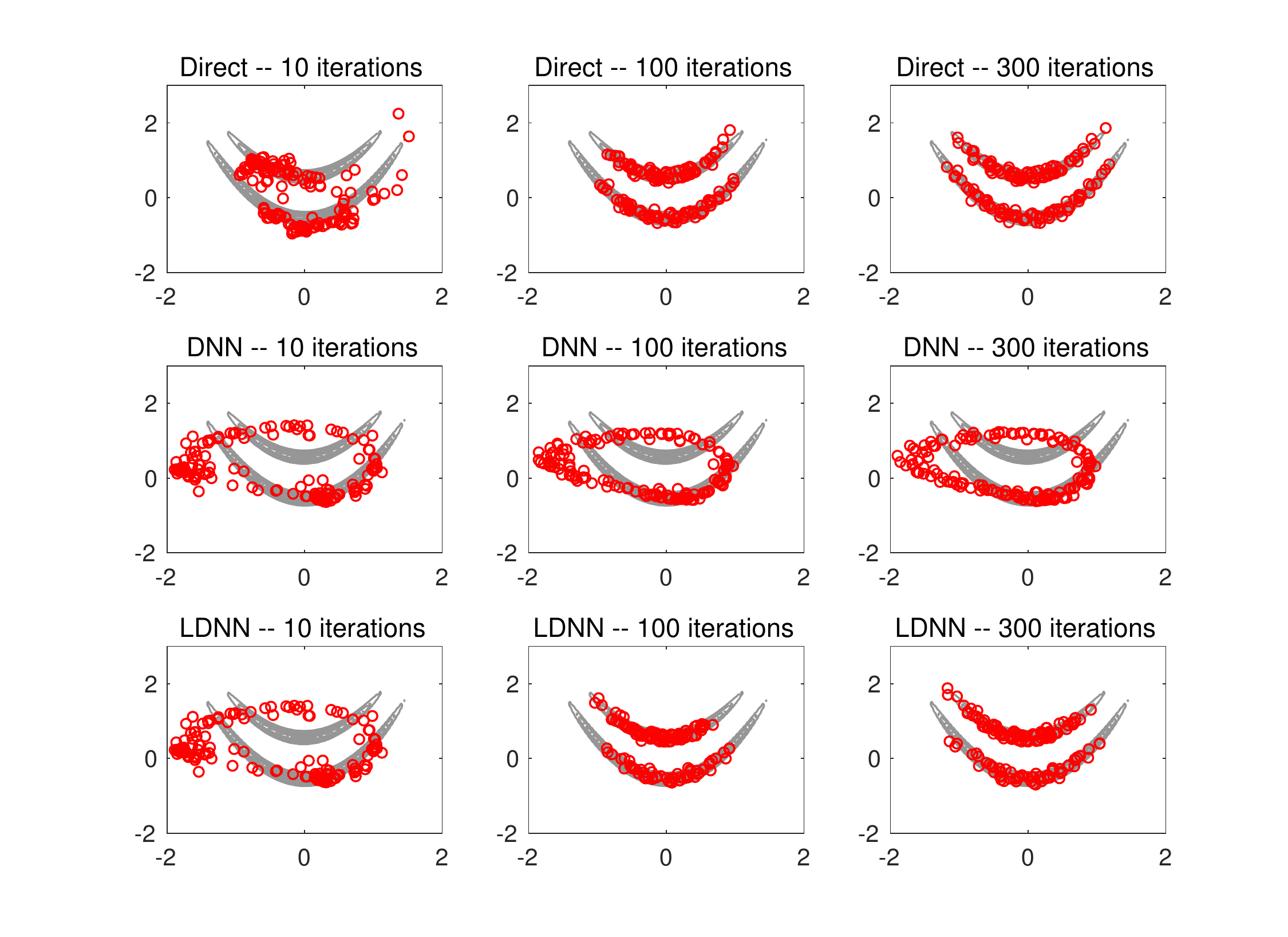} 
        \end{overpic}
     \end{center}
\caption{ The particles  obtained by various SVGD methods at the 10/100/300-th iterations on the double banana distribution. Here, the number of the iterations for LDNN is refer to the total number of iterations. 
}\label{resSVGD_eg1}
  \end{figure}

    \begin{figure}
\begin{center}
         \begin{overpic}[width=0.85\textwidth,trim= 35 10 45 15, clip=true,tics=10]{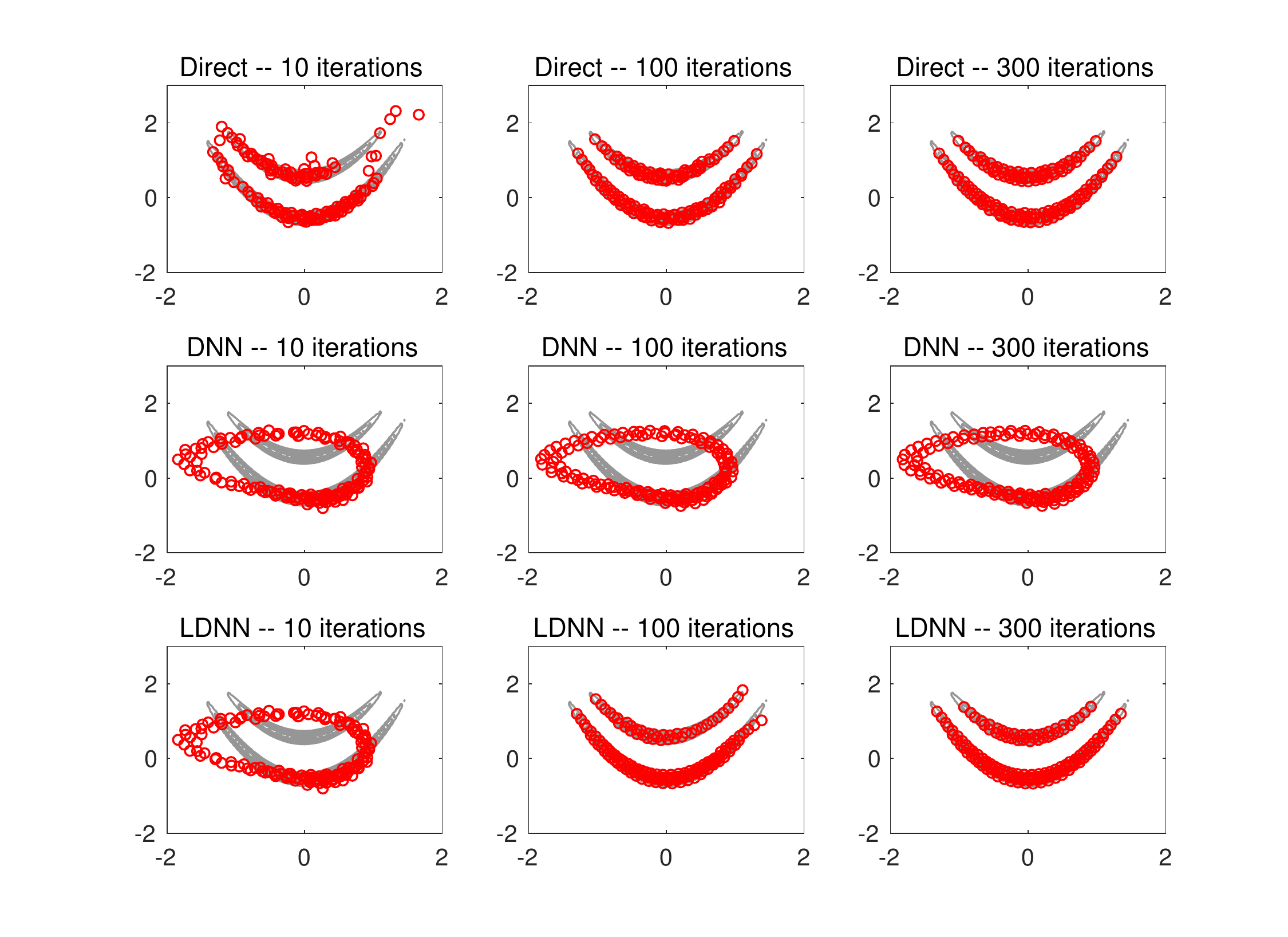} 
        \end{overpic}
     \end{center}
\caption{The particles  obtained by various SVN methods at the 10/100/300-th iterations on the double banana distribution. Here, the number of the iterations for LDNN is refer to the total number of iterations. 
}\label{resSVN_eg1}
  \end{figure}
  
Fig. \ref{resSVGD_eg1} provides the distribution of $N=100$ particles obtained by three different algorithms at selected iteration numbers, each with particles initially sampled from the two-dimensional standard Gaussian distribution.   The first row of  Fig. \ref{resSVGD_eg1} displays the performance of original SVGD, where the exact gradient  information is exploited  in the optimization.  At 100 iterations,  we can see that the particles  have already converged and the particles are all in the high-probability regions. The second row shows the performance of prior-based DNN algorithm, where the DNN is pre-trained in advance and kept unchanged during SVGD computations. Specially, we construct the DNN surrogate using 10 training points with 3 hidden layers and 20 neurons per layer. All the particles obtained by DNN quickly reach the high probability regions of the approximated posterior distribution, but they do not converge to the exact posterior.  Using this pre-trained DNN, we can refine it via Algorithm \ref{alg:LDNN}  and obtain the LDNN algorithm. The corresponding results are shown in the third row of Fig. \ref{resSVGD_eg1}. As expected, the particles spread over the support of the posterior  after 100 iterations, and all the particles reach the high probability regions of the exact target in the last iteration. Here, the number of the iterations for LDNN  refers to the total number of iterations, i.e.,  
$
\text{(No. of Outer Iterations)} \times \text{(No. of Inner Iterations).}
$

Our framework can easily to extend to other Stein variational inference algorithm. For example, we consider the Stein variational Newton method (SVN)\cite{Detommaso2018stein}, which generalize the SVGD algorithm by including second-order information. The corresponding results are shown in Fig. \ref{resSVN_eg1}.  It can be seen that  all the particles quickly reach the high probability regions of the exact or the approximated posterior distribution, due to the Newton acceleration in the optimization. Again, the  DNN algorithm gives a very poor density estimate, however, the LDNN algorithm yields a good approximation to the reference solution obtained by Direct scheme.  By comparing Figs. \ref{resSVGD_eg1} and \ref{resSVN_eg1}, we can obtain that the approximation results using LDNN are much more accurate than those of the prior-based DNN algorithm. 

\begin{figure}
\begin{center}
         \begin{overpic}[width=0.55\textwidth,trim= 35 10 45 15, clip=true,tics=10]{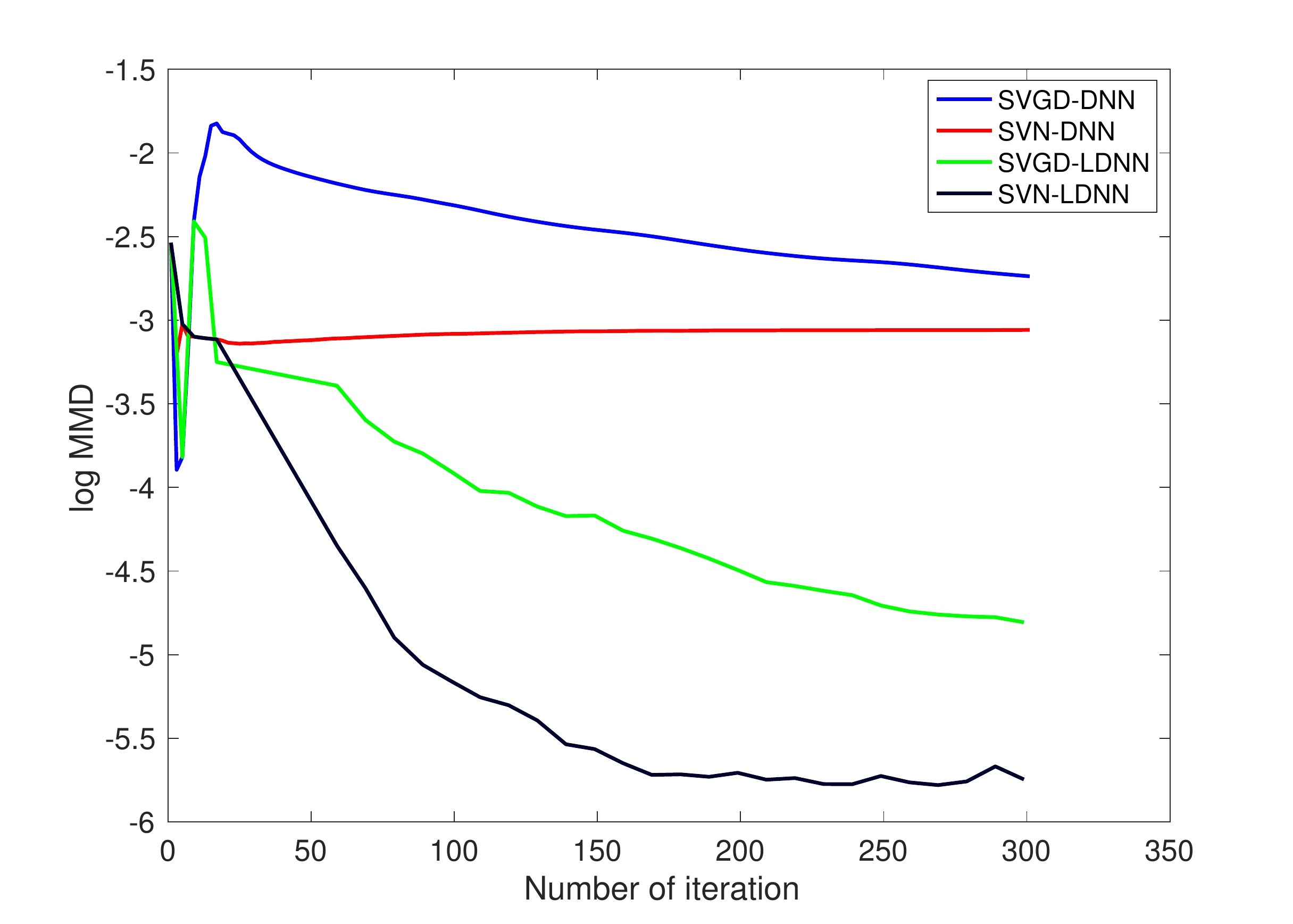} 
        \end{overpic}
     \end{center}
\caption{The log MMD vs. training iteration of different algorithms on the double banana distribution. }\label{reserr_eg1}
  \end{figure}

\begin{figure}
\begin{center}
         \begin{overpic}[width=0.85\textwidth,trim= 35 10 45 15, clip=true,tics=10]{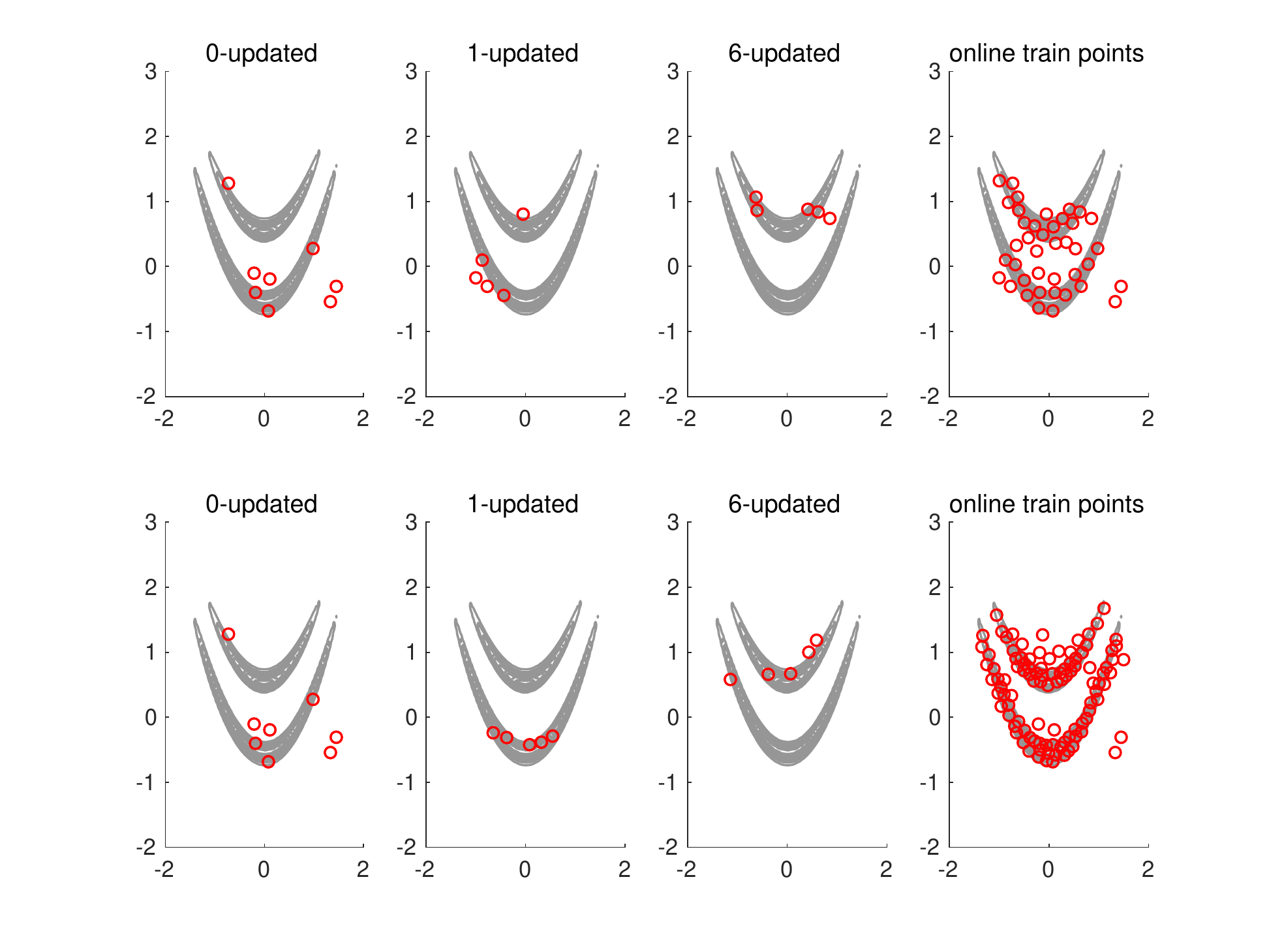}
        \put (9,9) {\scriptsize \bf SVN}
         \put (9,50) {\scriptsize \bf SVGD}
        \end{overpic}
     \end{center}
\caption{2D toy example. The evolution of the design pool by the LDNN algorithm. Red circles are new training points added to the training in each adaptation.  The right plot are the distribution of the total online training points. }\label{res_points_eg1}
  \end{figure}

To evaluate the sample quality, we compute the MMD  between the particles obtained by surrogate-based methods and the `true distribution' obtained by Direct SVN.  The MMD results for increasing iteration numbers are plotted in Fig. \ref{reserr_eg1}. It can be seen that LDNN offers improved performance over the DNN method. The main reason is the LDNN can gradually spread the online training points over the target density contour. The evolution of the design pool by the LDNN algorithm is shown in Fig. \ref{res_points_eg1}.

The computational costs, given by three-type mentioned algorithms, which is a count of the total number  of $N_{eval}$ times either the high-fidelity model $\mb{f}$ or $\nabla\mb{f}$ were evaluated, are presented in Table \ref{eg1_time}.   The main advantage of the DNN algorithm is that it does not require any high-fidelity model evaluations computed online.  For the LDNN, however, we indeed need the online high-fidelity model simulations. Nevertheless, in contrast to $150 \times 300$ high-fidelity gradient evaluations in the conventional SVGD (SVN), the number of high-fidelity model evaluations for the LDNN approach with SVGD (SVN) is only 50 (120).   As can be seen from the sixth and seventh column of Table \ref{eg1_time}, the LDNN approach does not significantly increase the computation time compared to the DNN approach, yet it offers significant improvement in the accuracy.  Despite these results being promising, it is important to note that updating the surrogate model may be costly itself. Note that the CPU times for the Direct approach are problem-dependent. The forward model in this example is chosen to be relatively inexpensive, and the online SVGD simulation is very cheap as its gradient has an analytical form.   For this reason, we envision that the proposed approach may have more potential when the computational solution of the forward problem is very expensive. 

   \begin{table}[tp]
      \caption{2D toy example. Computational times, in seconds, given by  different methods. }\label{eg1_time}
  \centering
  \fontsize{6}{12}\selectfont
  \begin{threeparttable}
    \begin{tabular}{ c cccccc}
  \toprule
 & \multicolumn{2}{c}{Offline}&\multicolumn{2}{c}{Online}\cr
\cmidrule(lr){2-3} \cmidrule(lr){4-5}

  \multirow{1}{*}{Method}  &$\text{$\#$ of model eval.}$&CPU(s) &$\text{$\#$ of model eval.}$&CPU(s)     &\multirow{1}{*}{Total time(s)}&\multirow{1}{*}{mmd}\cr
  \midrule
   SVGD                           & $-$       & $-$         & $150\times300 $          & 0.77      & 0.77  & $-$  \cr
   SVN                             & $-$         & $-$        &  $150\times300 $         &1.65         & 1.65  &$-$\cr
   SVGD-DNN                 & 10         & 2.6        &   $-$                               & 2.91        & 5.51    & 0.064\cr
   SVN-DNN                    & 10         & 2.6          &  $-$                             & 3.71       & 6.31    & 0.047\cr
   SVGD-LDNN              & 10         & 2.6          &  50                              & 8.46         & 11.06   & 0.0082\cr
   SVN-LDNN                 & 10         & 2.6          & 120                               &32.25       &  34.85   &0.0033\cr
      \bottomrule
      \end{tabular}
    \end{threeparttable}

\end{table}
  
  \begin{figure}
\begin{center}
         \begin{overpic}[width=0.45\textwidth,trim= 25 0 45 15, clip=true,tics=10]{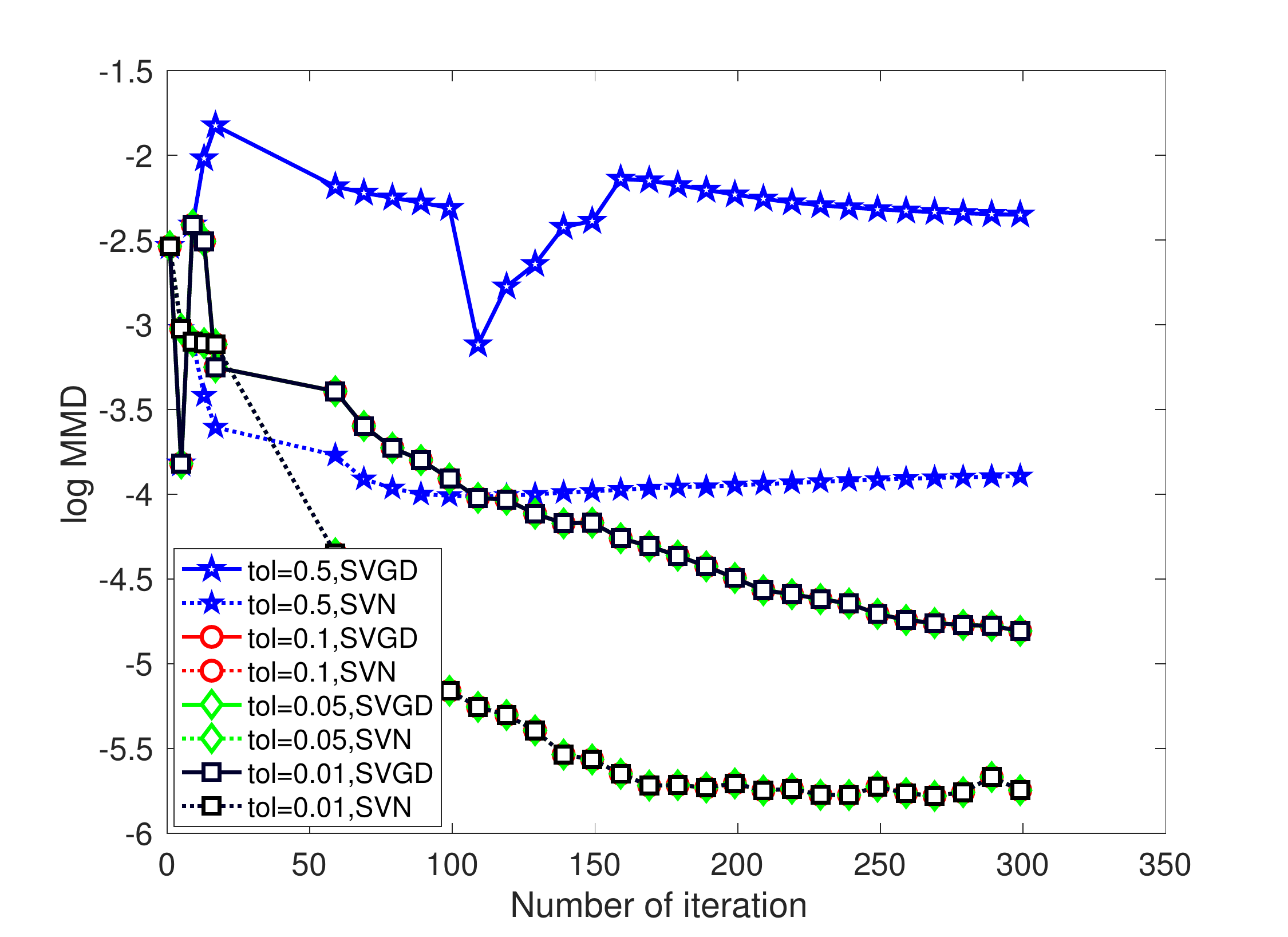}
        \end{overpic}
                \begin{overpic}[width=0.45\textwidth,trim= 20 0 30 15, clip=true,tics=10]{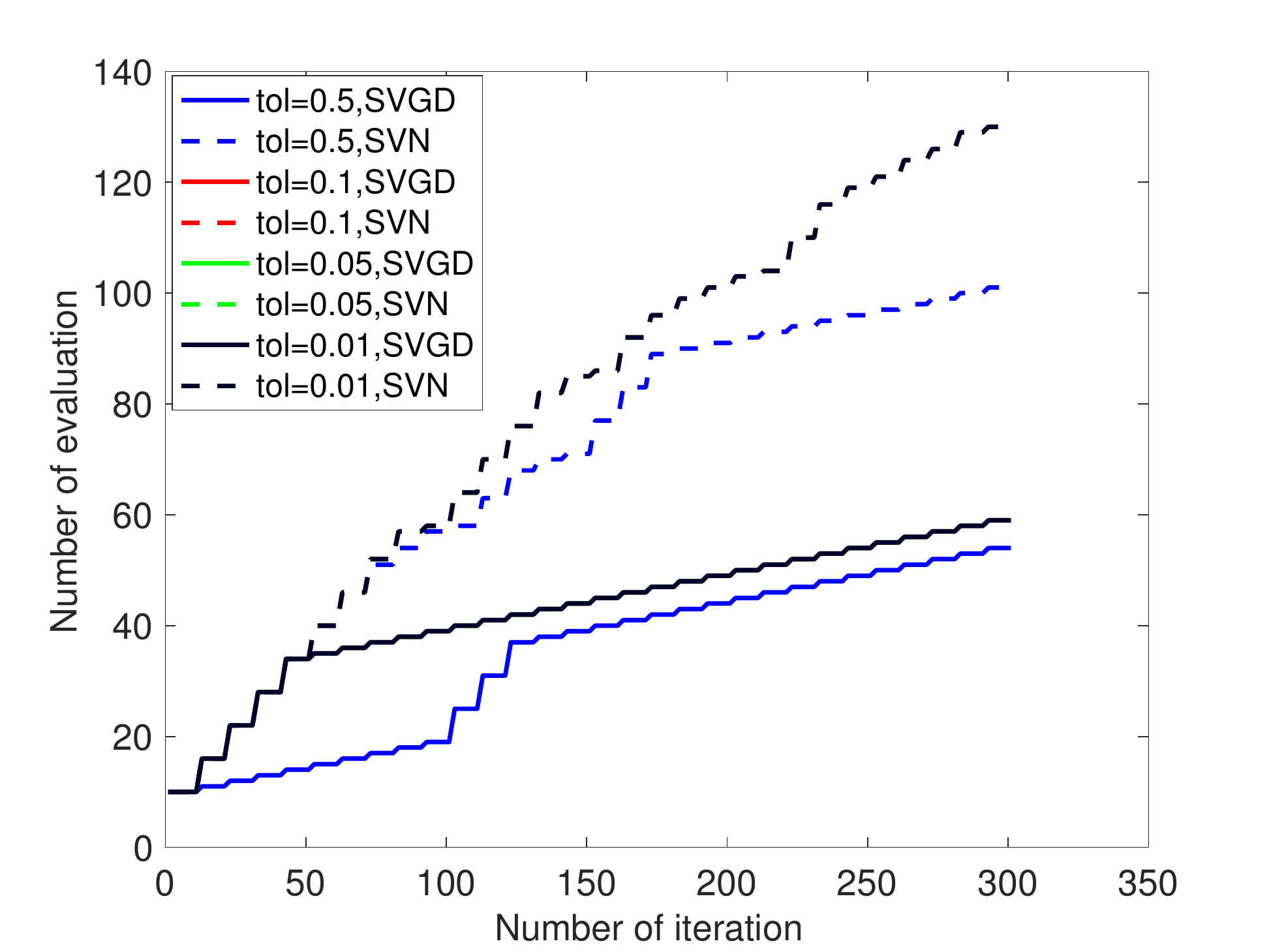}
        \end{overpic}
     \end{center}
\caption{ The log-MMD and the total number of high-fidelity model evaluations of LDNN vs. training iteration of different algorithms using various numbers of threshold tol.  }\label{resLDNN_eg1_tol}
  \end{figure}
  
    \begin{figure}
\begin{center}
         \begin{overpic}[width=0.45\textwidth,trim= 35 0 45 15, clip=true,tics=10]{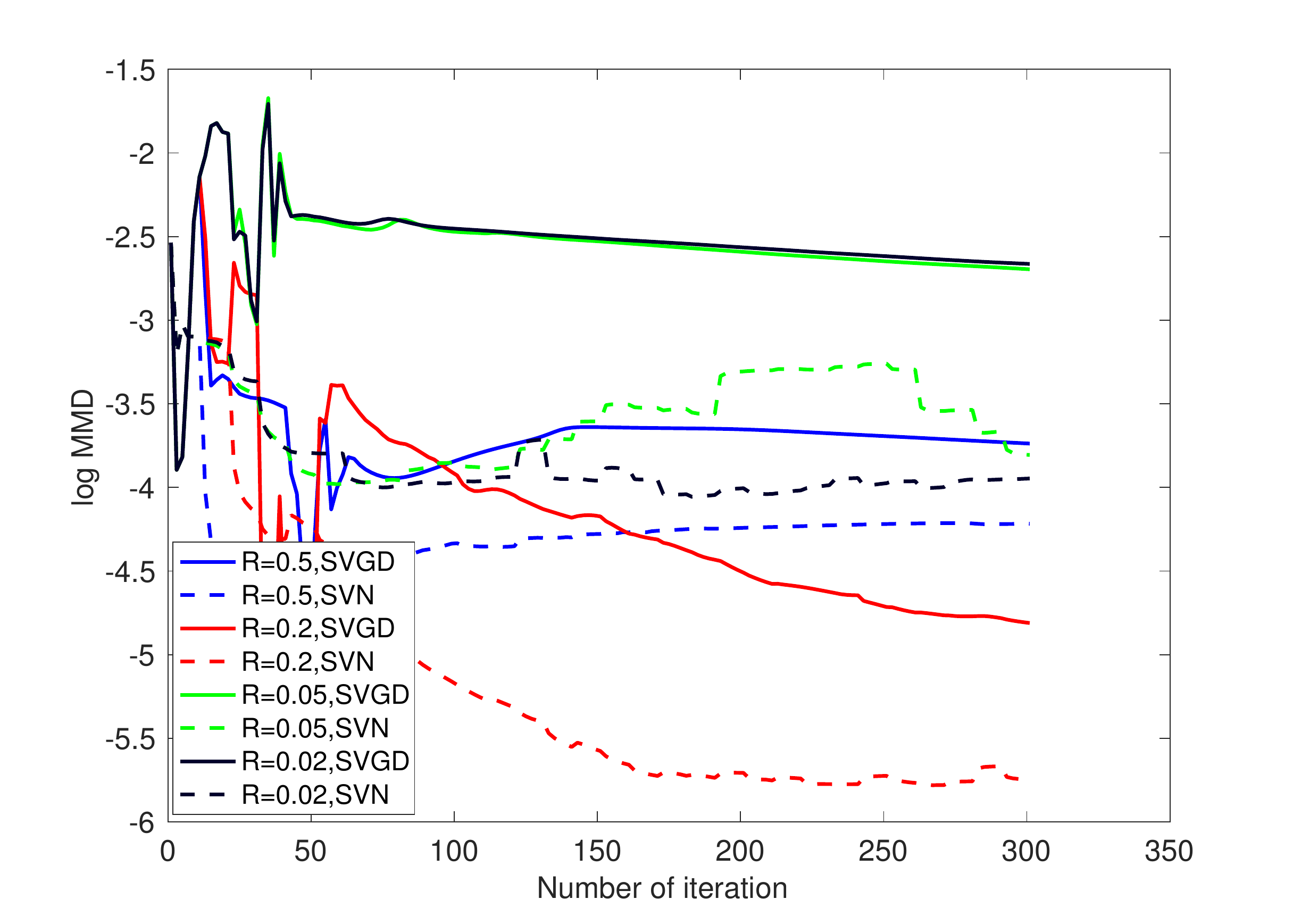}
        \end{overpic}
                \begin{overpic}[width=0.45\textwidth,trim= 35 0 45 15, clip=true,tics=10]{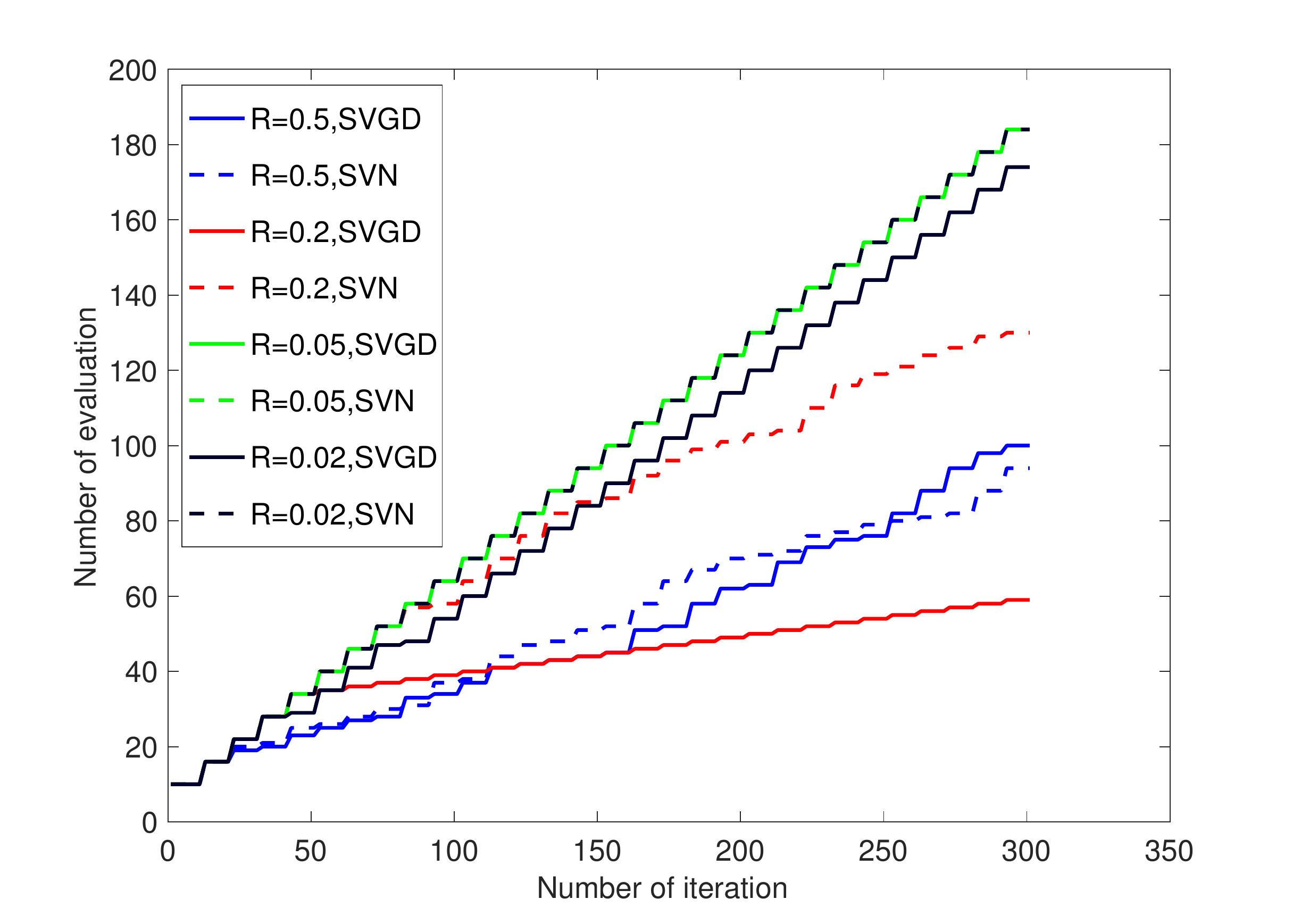}
        \end{overpic}
     \end{center}
\caption{ The log-MMD and the total number of high-fidelity model evaluations of LDNN vs. training iteration of different algorithms using various numbers of parameter $R$.  }\label{resLDNN_eg1_R}
  \end{figure}
  
  \begin{table}[h]
  \centering
 % \fontsize{8}{10}\selectfont
\begin{tabular}{|l|l|l|l|l|l|}
\hline
\diagbox{L}{$d_k$} & 10& 20 & 30 & 50 &80\\
\hline
1 & 0.0219 (107) & 0.0029 (70) & 0.0063 (54) & 0.0029 (93)  & 0.0027 (73) \\
\hline
2 & 0.0215 (105) & 0.0012 (75) & 0.0059 (71) & 0.0046 (70)  & 0.0106 (109)\\
\hline
3 & 0.0072 (97) & 0.0031 (96) & 0.0069 (95) & 0.0053 (72) & 0.0249 (105) \\
\hline
\end{tabular}
  \caption{Example 1:  The MMD (the number of online high-fidelity model evaluations) obtained using  LDNN-based SVN approach with  $Q=5$.}\label{rel_L}
    \centering
 % \fontsize{8}{10}\selectfont
\begin{tabular}{|l|l|l|l|l|l|}
\hline
\diagbox{Q}{$d_k$} & 10& 20 & 30 & 50& 80 \\
\hline
5 &  0.0072 (97) & 0.0031 (96) & 0.0069 (95) & 0.0053 (72) & 0.0249 (105)\\
\hline
10 & 0.0040 (68) & 0.0045 (91) & 0.0066 (98) & 0.0052 (120)  & 0.0212 (96)\\
\hline
15 & 0.0050 (73) & 0.0084 (105) & 0.0064 (118) & 0.0075 (87)  & 0.0079 (100) \\
\hline
20 & 0.0014 (50) & 0.0086 (118) & 0.0061 (116) & 0.0043 (114)  & 0.0016 (107) \\
\hline
\end{tabular}
  \caption{2D toy example:  The MMD (the number of online high-fidelity model evaluations) obtained using LDNN-based SVN approach with $L=3$. }\label{rel_Q}
\end{table}
  
Next,  we investigate the influence of the tuning parameters $tol$ and $R$.  Intuitively, one would expect that the accuracy of the LDNN will improve as the value of the threshold $tol$ decreases. To verify this proposition, we test several constant values choosing from $tol \in \{0.5,0.1,0.05,0.01\}$. The MMD and the total number of high-fidelity model evaluations of LDNN with respect to training iterations  using various numbers of threshold $tol$ are displayed in Fig.\ref{resLDNN_eg1_tol}. It can be seen that the MMD decreases when threshold $tol$ gets smaller; these values trigger more frequent refinements.  Interestingly, the numerical results are practically the same for $tol<0.1$ in the case of Example 1.  This indicates that accurate numerical results can be  obtained when using even a relatively large number of $tol$. The sensitivity of the numerical results with respect to the parameter $R$ is shown in Fig. \ref{resLDNN_eg1_R}.  We can see that the accuracy of the numerical results improves as the radius $R$ decreases.  However, if  $R$ is too small, the results becomes  inaccuracy as it  introduces tightly clustered points.   To reduce the online computational cost and retain the accuracy of estimate results, a reasonable choice of the size of $R$ may be $0.05<R<0.5$. 

In order to analyze the sensitivity of the proposed method with respect to the structures of the DNN, we first  consider a prior-based DNN surrogate  with $L\in\{1,2,3\}$ hidden layers and $d_k \in \{10,20,30,50,80\}$ neurons per layer using $N=20$ training points.  Using those pre-trained DNN, we can run our LDNN approach. The MMD (and the required number of the online high-fidelity model evaluations) of the LDNN-based SVN approach for Example 1 are presented in Table \ref{rel_L}. As shown in this table, the computational results for LDNN with different depth $L$ and width $d_k\in\{20,30,50\}$ are almost the same.  The numerical results using various values of $Q$  are shown in Table \ref{rel_Q}.  As expected, the LDNN approach admits reasonably accurate results for different values of $Q$.  These numerical results  also demonstrate the robustness of the LDNN approach.

\subsection{2D heat source inversion}\label{heatsource}
Next our attention turns to whether LDNN-based SVGD improves over the original SVGD method and how it compares to the DNN-based SVGD when the system is modeled by PDEs. To this end, we consider a 2D heat source inversion problem which adapted from \cite{yan+guo2015}. 
Consider the following model in the domain $D=[0,1]\times [0,1]$
\begin{eqnarray}\label{2dsource}
\begin{array}{rl}
^cD_t^{\alpha}u(\textbf{s},t)-\triangle u(\textbf{s},t)&=e^{-t}\exp\Big[-0.5\Big(\frac{\|x-\mathbf{s}\|}{0.1}\Big)^2\Big],\quad D\times [0, 1],\\ 
\nabla u \cdot \textbf{n}&=0, \quad \mbox {on} \,\partial{D},\\
u(\textbf{s},0)&=0, \quad  \mbox{in}\, D.
 \end{array}
\end{eqnarray}
The goal is to determine the source location $x=(x_1,x_2)$ from noisy measurements of the $u$-field at a finite set of locations and times.  Here $^c D^{\alpha}_t \,(0<\alpha<1)$ denotes the Caputo fractional derivative of order $\alpha$ with respect to $t$  and it is defined by 
\begin{eqnarray*}
^cD_t^{\alpha}u(\mathbf{s},t)=\frac{1}{\Gamma(1-\alpha)} \int ^t_0 \frac{\partial{u(\mathbf{s},\eta)}}{\partial{\eta}} \frac{d\eta}{(t-\eta)^{\alpha}}, \, 0<\alpha<1,
\end{eqnarray*}
where $\Gamma (\cdot)$ is the Gamma function.  Note that the gradient  of the target distribution in this example  cannot be easily calculated analytically, and at the same time the  evaluation of the froward model are computationally expensive too. In such a setting, the standard SVGD  algorithm cannot be applied.

\begin{figure}
\begin{center}
     \begin{overpic}[height=4.4cm,width=3.3cm,trim= 35 0 45 15, clip=true,tics=10]{figures/initial_points2-eps-converted-to.pdf}
           \put (4,99) {\scriptsize {\bf Initial design  points}}
          \put (43,15) {\footnotesize \bf Case 1}
        \end{overpic}
         \begin{overpic}[height=4.4cm,width=3.3cm,trim= 35 0 45 15, clip=true,tics=10]{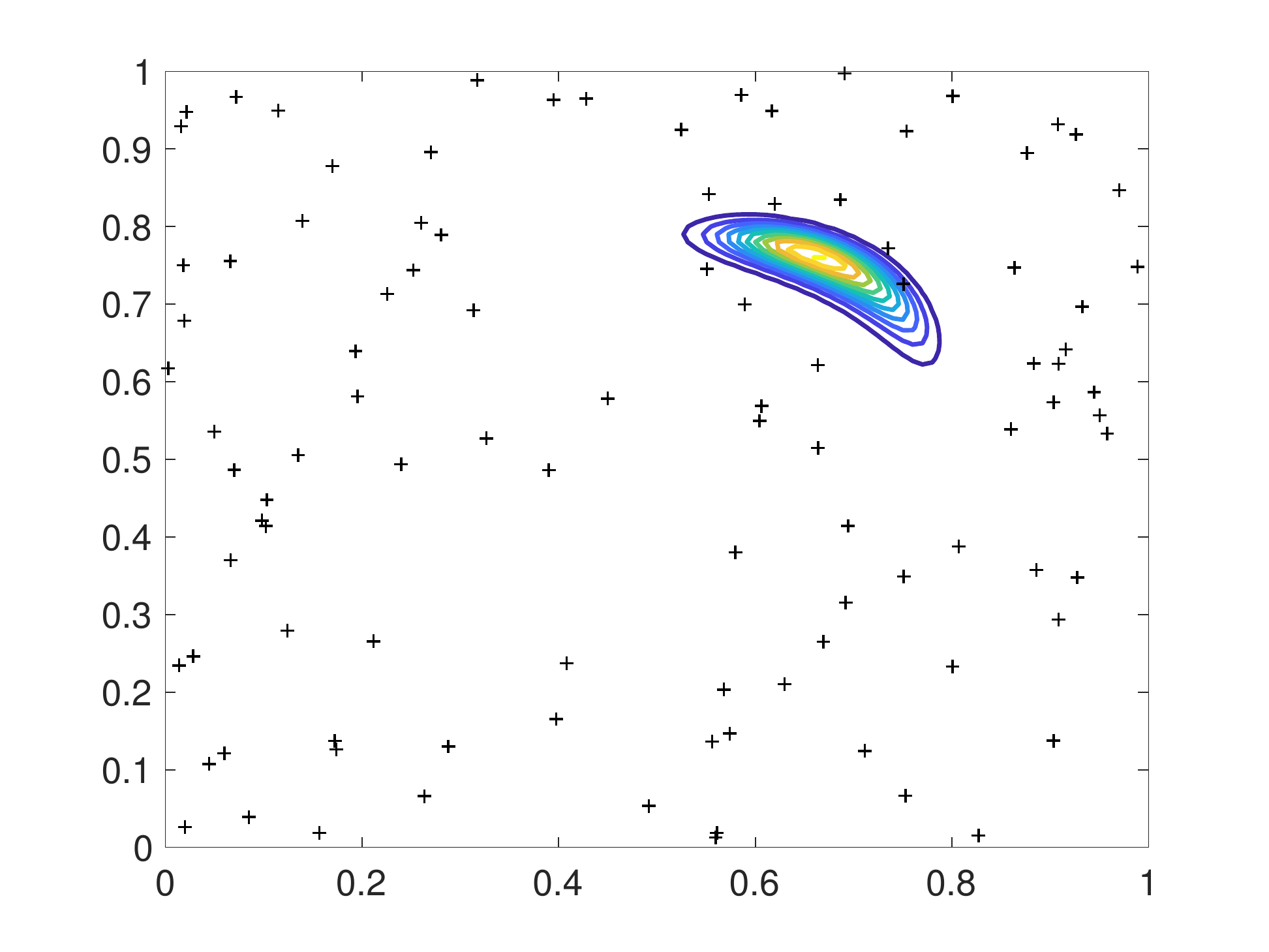}
           \put (14,99) {\scriptsize {\bf Initial particles}}
        \end{overpic}
    \begin{overpic}[height=4.4cm,width=3.3cm,trim= 35 0 45 15, clip=true,tics=10]{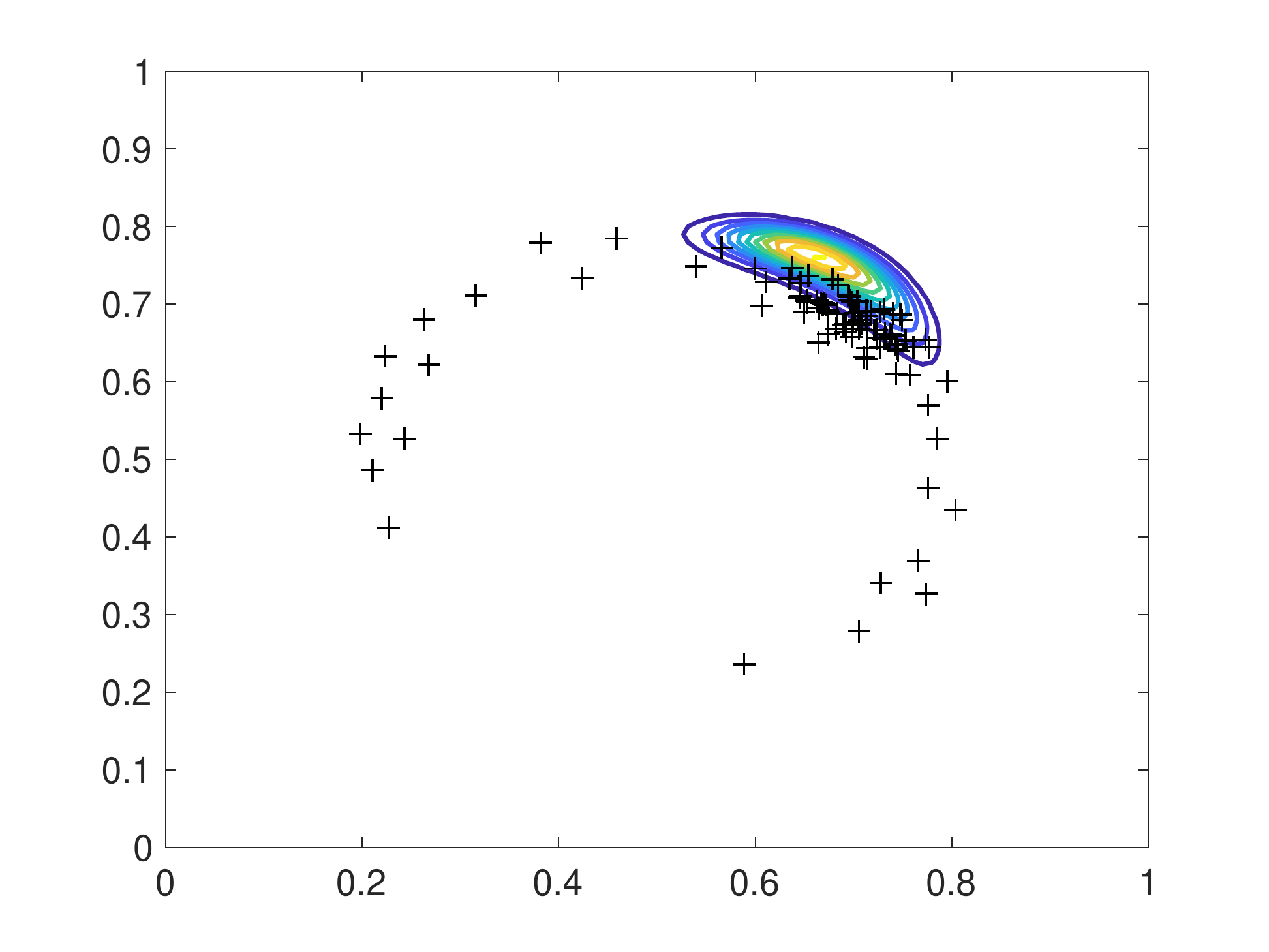}
       \put (4,99) {\scriptsize {\bf Final particles-DNN}}
  \end{overpic}
      \begin{overpic}[height=4.4cm,width=3.3cm,trim= 35 0 45 15, clip=true,tics=10]{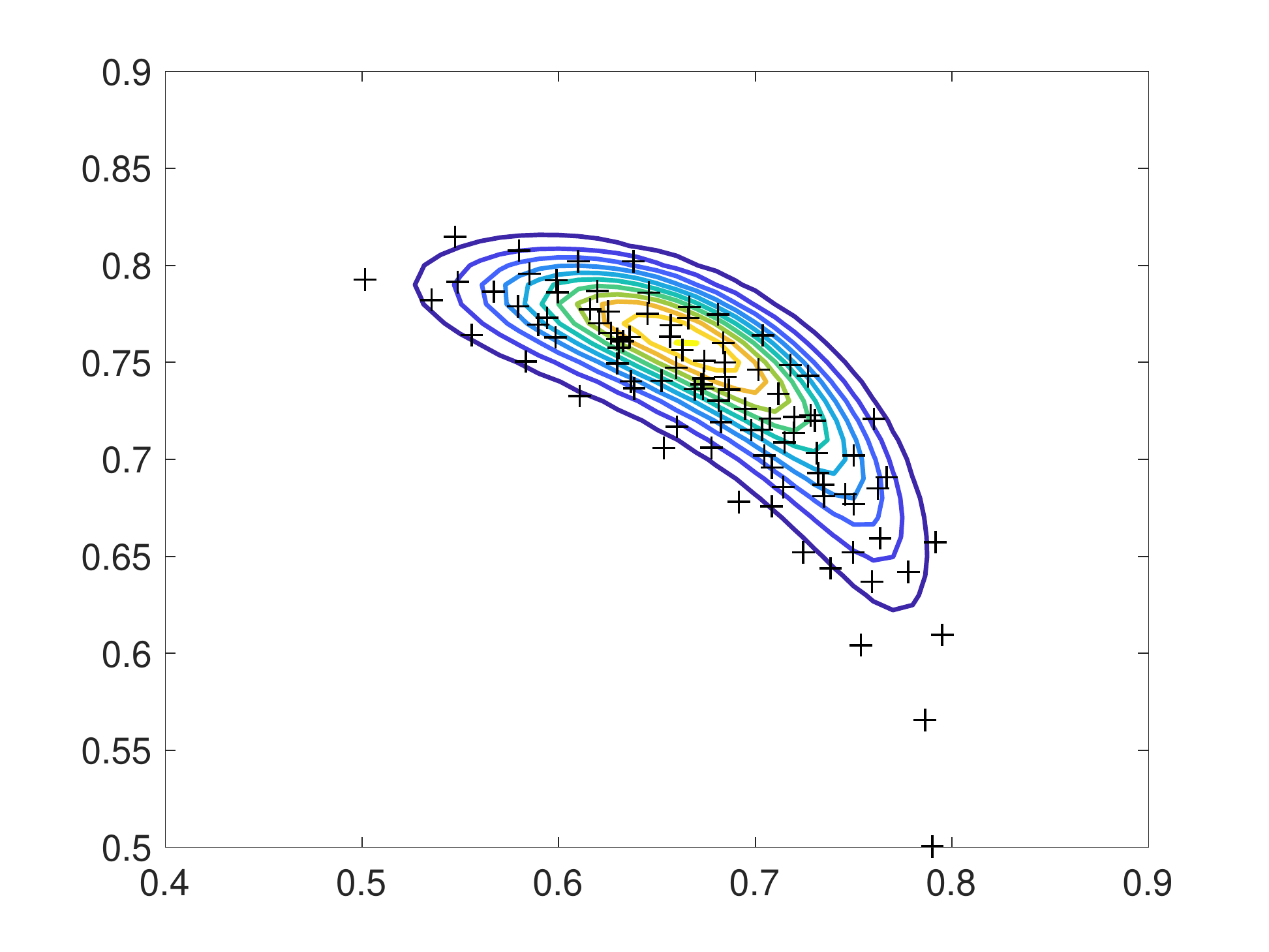}
        \put (4,99) {\scriptsize {\bf Final particles-LDNN}}
  \end{overpic}
     \begin{overpic}[height=4.4cm,width=3.3cm,trim= 35 0 45 15, clip=true,tics=10]{figures/initial_points1-eps-converted-to.pdf}
         \put (43,15) {\footnotesize \bf Case 2}
  \end{overpic}
      \begin{overpic}[height=4.4cm,width=3.3cm,trim= 35 0 45 15, clip=true,tics=10]{figures/initial_partcle-eps-converted-to.pdf}
  \end{overpic}
    \begin{overpic}[height=4.4cm,width=3.3cm,trim= 35 0 45 15, clip=true,tics=10]{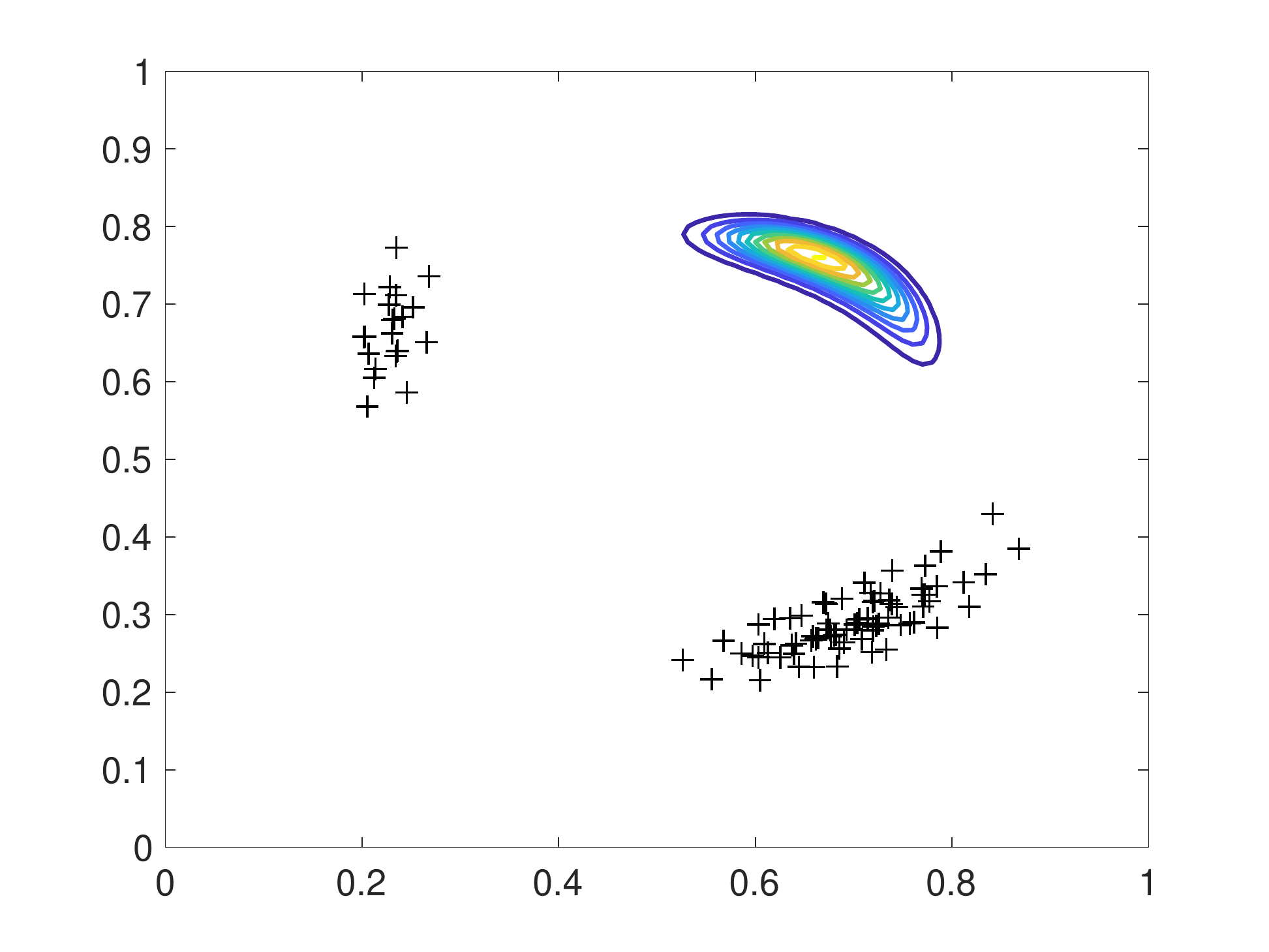}
  \end{overpic}
      \begin{overpic}[height=4.4cm,width=3.3cm,trim= 35 0 45 15, clip=true,tics=10]{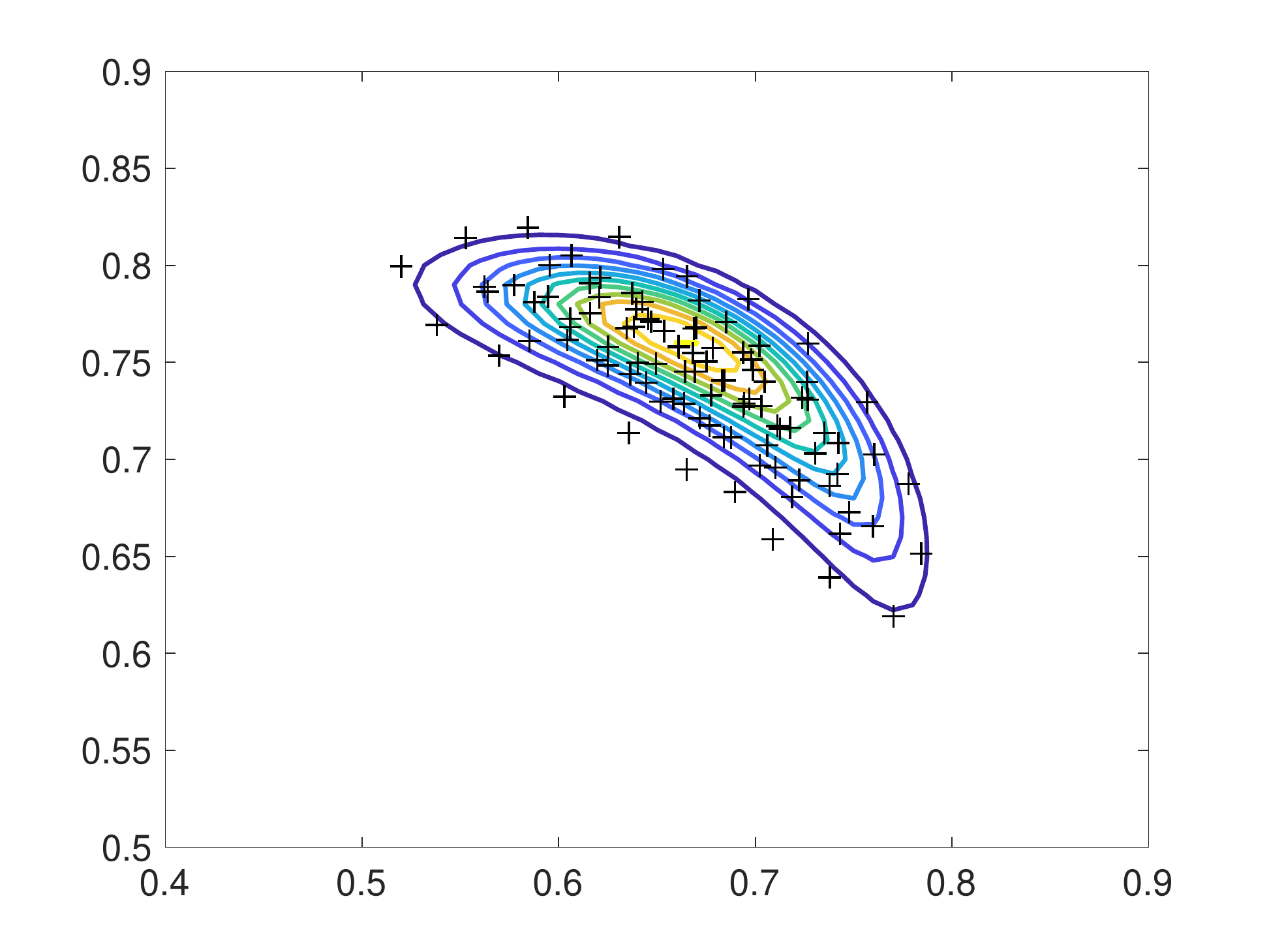}
  \end{overpic}
  \end{center}
\caption{The results of 2D heat source inversion. }\label{res_eg2}
  \end{figure}
  
\begin{figure}
\begin{center}
         \begin{overpic}[width=0.45\textwidth,trim= 35 0 45 15, clip=true,tics=10]{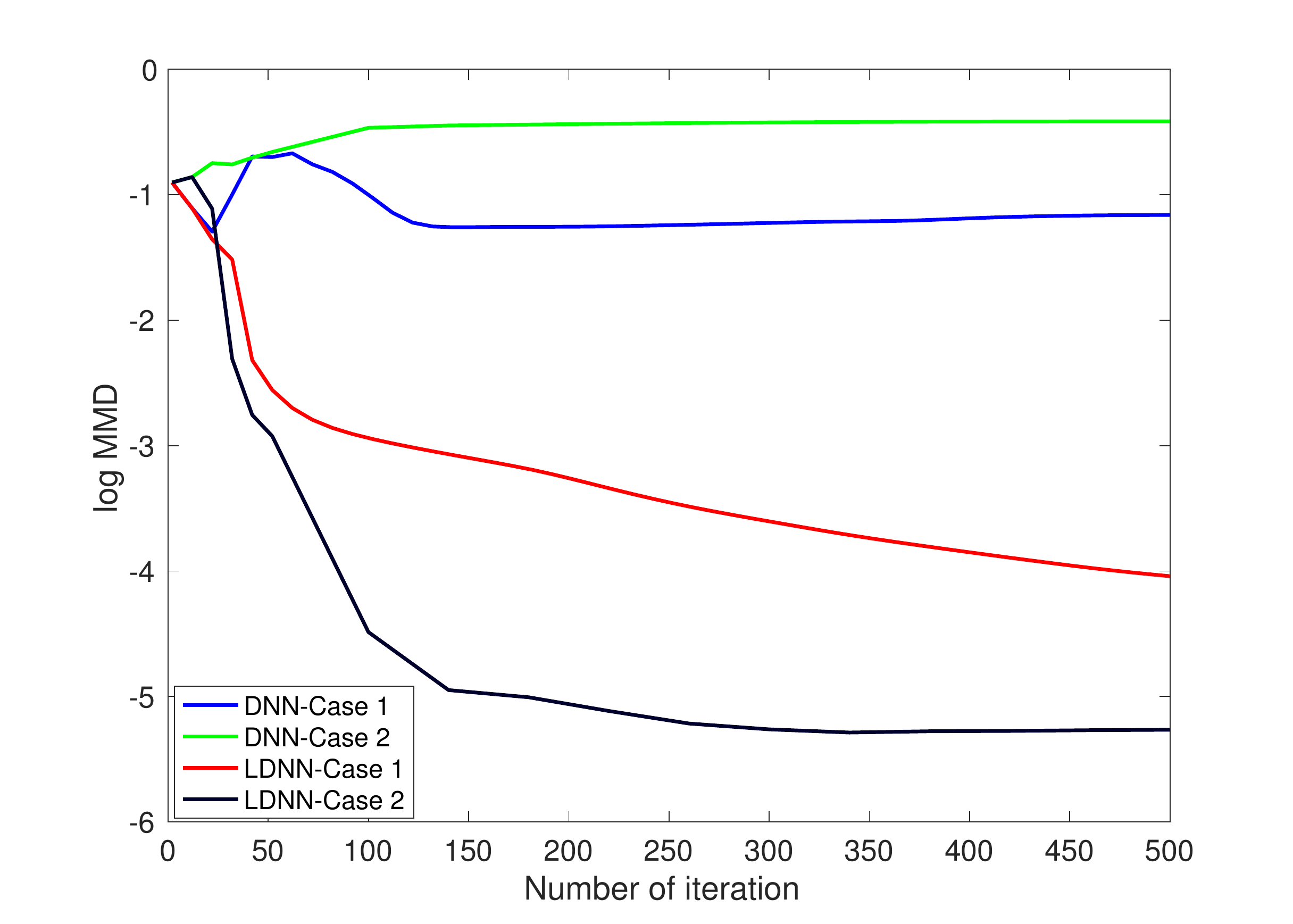}
        \end{overpic}
                \begin{overpic}[width=0.43\textwidth,trim= 30 0 35 15, clip=true,tics=10]{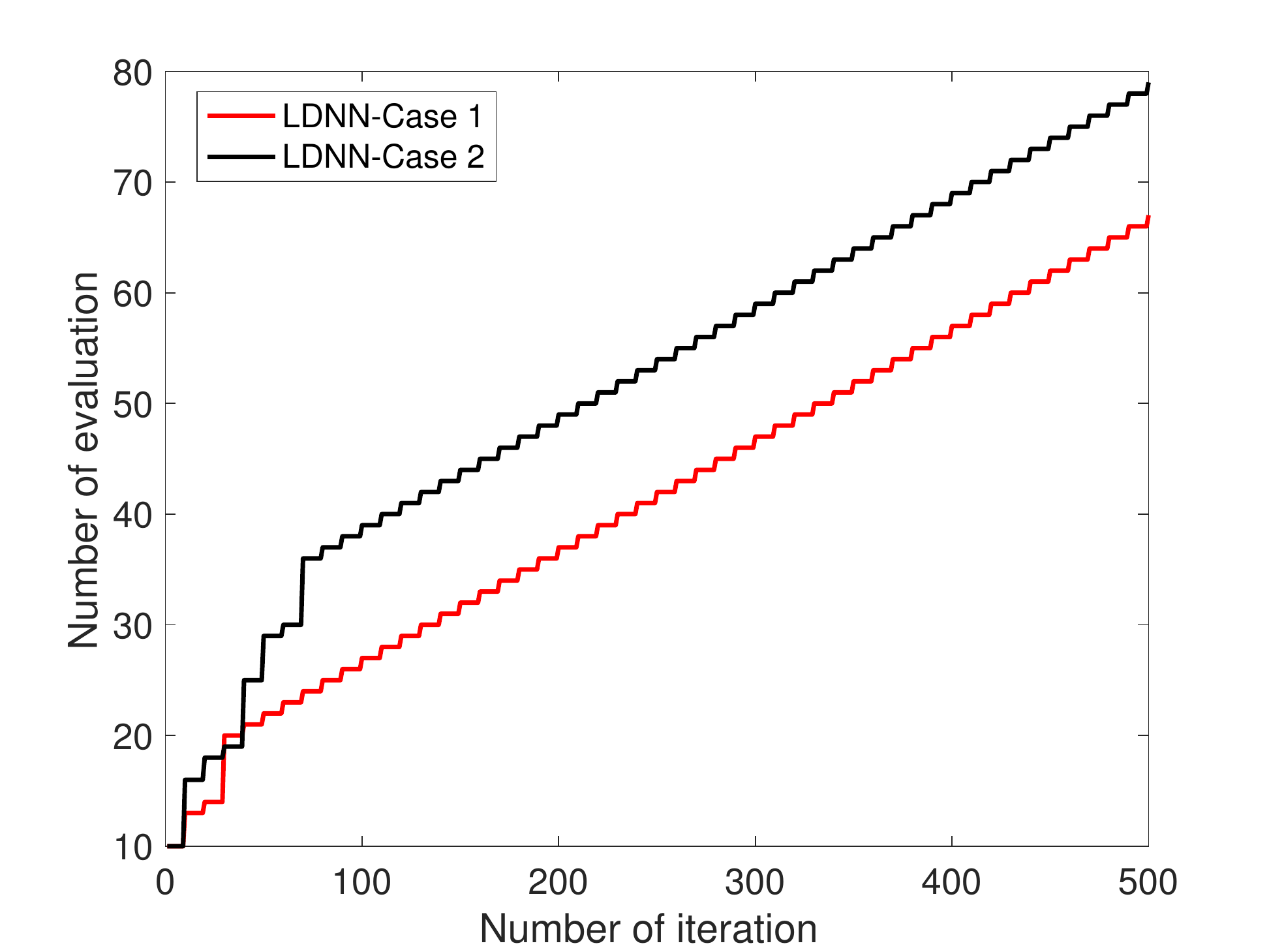}
        \end{overpic}
     \end{center}
\caption{ The log-MMD and the total number of high-fidelity model evaluations of LDNN vs. training iteration.}\label{resLDNN_eg2}
  \end{figure}
  
%  \begin{figure}
%\begin{center}
%         \begin{overpic}[width=0.5\textwidth,trim= 10 0 45 15, clip=true,tics=10]{figures/err_x_eg2-eps-converted-to.pdf}
%        \end{overpic}
%     \end{center}
%\caption{The relative errors with number of iteration for LDNN approach.  }\label{resLDNN_eg2_x.}
%  \end{figure}

To obtain the high-fidelity data, we solve the equation (\ref{2dsource}) using a finite difference/ spectral approximations (\cite{Lin+Xu2007}) with time step $\Delta t=0.01$ and polynomial degree $P=6$.  Noisy point-wise observations are taken from the solution of the field of the PDE at a uniform $3 \times 3$ sensor network. At each sensor location, two measurements are taken at time $t = 0.25$ and $t = 0.75$, which corresponds to a total of  18 measurements.  The prior on $x$ is set as $x_i\sim U(0,1)$.  The likelihood assumes additive and i.i.d. Gaussian errors  with mean zero and variance 0.04 for each observation. 
In order not to commit an `inverse crime',  we generate the data by solving the forward problem at a much higher resolution than that used in the inversion, i.e., with  $P=10$. In the examples below, unless otherwise specified, the Caputo fractional derivative of order $\alpha$ is 0.5.   To evaluate how well the particles approximate the posterior, we first obtain a ground truth set of $10^4$ MCMC samples from a long run of an efficient adaptive MCMC\cite{Haario+Laine+Mira2006}.  Given this reference set of samples, the MMD between the particles and the MCMC samples was used to assess closeness of all empirical measures to the target $\pi(x)$. 

Similar to the first example, we numerically investigate the efficiency of the LDNN approach. Using the same setting as Example 1, we  consider two type initial design points to train the prior-based DNN: (1)  starting with $n_t=10$ nodes randomly chosen in $[0,1] \times [0,1]$; (2)  starting with $n_t=10$ nodes randomly initialization in $[0, 0.5] \times [0, 0.5]$ which is far away from the exact target, see the first column of Fig.\ref{res_eg2}.  In the numerical experiment, we use $N=100$ particles and the same initial particles sampled from the the uniform prior distribution. The online training points choosing by LDNN are shown in Fig.\ref{uppoints_eg2}.  The final results in Fig.\ref{res_eg2} show that the particles returned by LDNN approximate the target distributions reasonably well, even staring with a bad design.  Again, the prior-DNN approach can not obtain  accurate results. 

The left figure of Fig. \ref{resLDNN_eg2} shows the decay of the MMD  for the DNN approximation constructed at the initial step of SVGD, in contrast with the LDNN approach. The total number of high-fidelity model evaluations of LDNN approach are shown in right figure of Fig. \ref{resLDNN_eg2}.  It is clearly seen that the prior-based DNN approach admits very large approximation error, especially for Case 2,  due to the fact that the exact target is far away from what is assumed in the location of the initial design pool.   This demonstrates the advantage of the online adaptive refinement in terms of accuracy of the LDNN approach.

We report the computational cost of DNN and LDNN  approach in the SVGD process up to step $N_{iter} =500$ in Table \ref{eg2_time}.   The more challenging nature of this experiment meant accurate computation of the original SVGD was precluded, due to the fact that a sufficiently high-quality empirical approximation of $\nabla \log \pi(x)$ could not be obtained. For ease of comparison, we use the cost of  $N_{iter}\times N$ high-fidelity model evaluations to represent the CPU time of the original SVGD approach.  In this sense, the CPU time of evaluating the conventional SVGD is more than 1210s, while the CPU time of prior-based DNN method for Case 1 is about 6.84s.  Although the a prior-based DNN approach can gain the computational efficiency, the estimation accuracy cannot be guaranteed.  In contrast, for the LDNN approach, the offline and online CPU times are 3.33s and 17.73s, respectively, meaning that the LDNN approach can provide  much more accurate results, yet with less computational time.  From the results we can see that  the LDNN achieves mostly over 50X speedup compared to the Direct approach  in terms of CPU time.  Moreover, with larger number of  particles, the LDNN leads to similar number of online heigh-fidelity model evaluations for the same  tolerance $tol$, and can achieve higher speedup since the online time to refine the  DNN  does not change much.  This confirms the efficiency of the LDNN approach.

  \begin{table}[tp]
      \caption{2D heat source inversion. Computational times given by three different methods. }\label{eg2_time}
  \centering
  \fontsize{6}{12}\selectfont
  \begin{threeparttable}
    \begin{tabular}{l cccccc}
  \toprule
 & \multicolumn{2}{c}{Offline}&\multicolumn{2}{c}{Online}\cr
\cmidrule(lr){2-3} \cmidrule(lr){4-5}

  \multirow{1}{*}{Method}  &$\text{$\#$ of model eval.}$&CPU(s) &$\text{$\#$ of model eval.}$&CPU(s)     &\multirow{1}{*}{Total time(s)}&\multirow{1}{*}{mmd}\cr
  \midrule
   SVGD                           & $-$       & $-$         & $100\times500 $          & $\sim$1210      & $\sim$1210  & $-$  \cr
 
 Case 1-DNN                  & 10         & 3.33        &   $-$                           & 3.51        & 6.84    & 0.304\cr
 Case 2-DNN                 & 10         & 3.23         &  $-$                           & 3.51         & 6.74    & 0.607\cr
Case 1-LDNN                & 10         & 3.33          &  72                            & 17.73         & 21.06    & 0.0182\cr
 Case 2-LDNN              & 10         & 3.23           &  66                            &15.32         & 18.55    &0.0043\cr
      \bottomrule
      \end{tabular}
    \end{threeparttable}

\end{table} 

 \subsection{Estimating the diffusion coefficient} 
\begin{figure}
\begin{center}
         \begin{overpic}[width=0.5\textwidth,trim= 10 0 45 15, clip=true,tics=10]{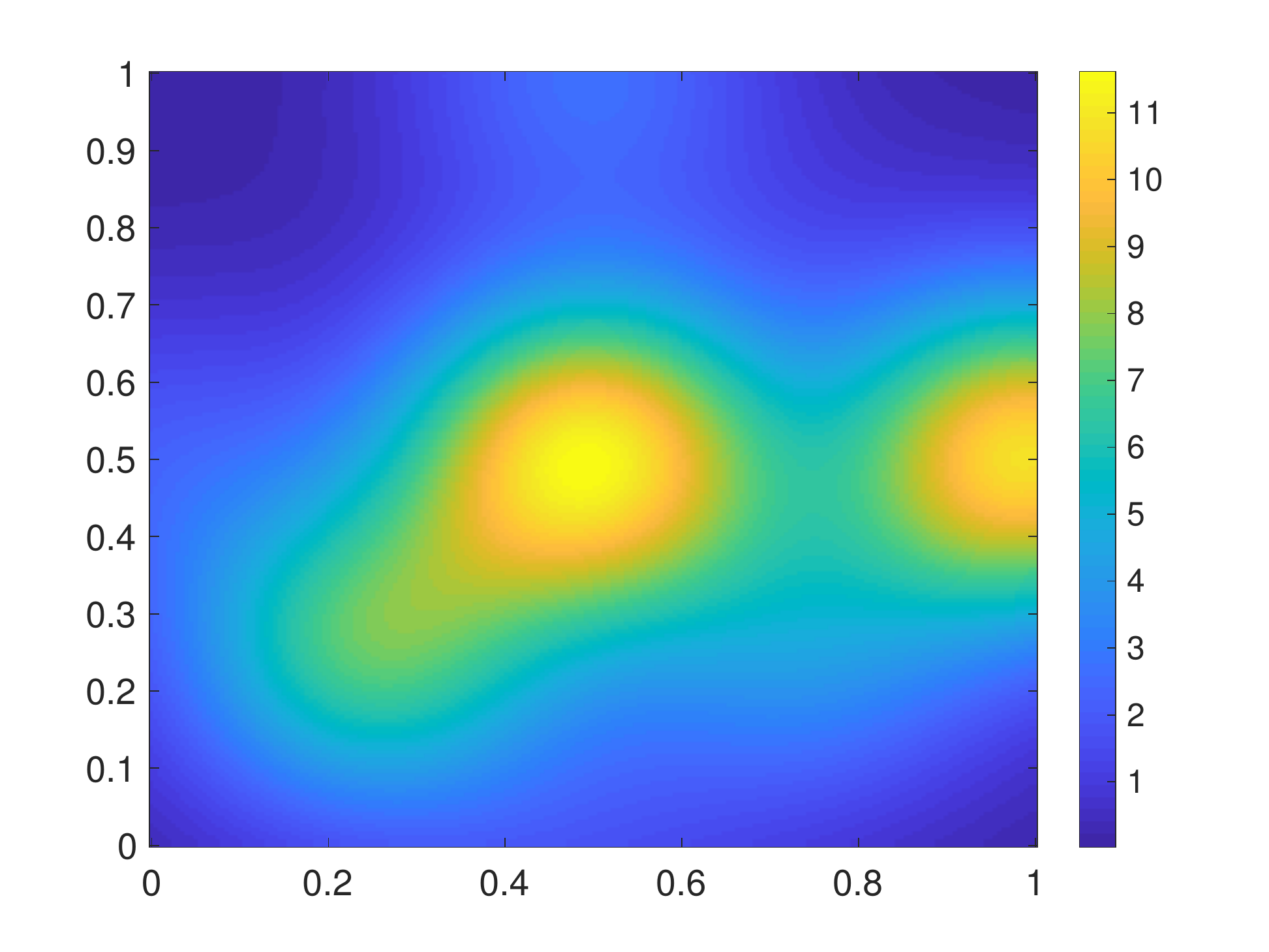}
        \end{overpic}
     \end{center}
\caption{ The true permeability used for generating the synthetic data set.  }\label{exact_eg3}
  \end{figure}

In the last example, we illustrate the LDNN approach on the nonlinear inverse problem of estimating the  diffusion coefficient. Consider the following two dimensional time-fractional PDEs \cite{Yan+Zhou2019PCEKI}
\begin{eqnarray}\label{2dtfpde}
\begin{array}{rl}
^cD_t^{\alpha}u-\nabla\cdot(\kappa(\mb{s}) \nabla u(\mb{s},t))&=e^{-t}\exp\Big(-\frac{\|\mb{s}-(0.25,0.75)\|^2}{2\times0.1^2}\Big),\quad \Omega\times [0, 1],\\
\nabla u \cdot \textbf{n}&=0, \quad \mbox {on} \,\partial{\Omega},\\
u(\mb{s},0)&=0, \quad  \mbox{in}\, \Omega.
 \end{array}
\end{eqnarray}
The goal is to determine the diffusion coefficient  $\kappa(\mb{s})$ from noisy measurements of the $u$-field at a finite set of locations and times.
We consider the following permeability field $\kappa(\mb{s};x)$
\begin{eqnarray*}
\kappa(\mb{s}; x)=\sum^{9}_{i=1}x_i \exp(- 0.5\frac{\|\mb{s}-\mb{s}_{0,i}\|^2}{0.15^2}),
\end{eqnarray*}
where $\{\mb{s}_{0,i}\}^{9}_{i=1}$ are the centers of the radial basis function.  The prior distributions on each of the weights $x_i, i=1,\cdots, 9$ are independent and log-normal; that is, $\log(x_i)\sim N(0,1)$.   The true permeability field used to generate the test data is shown in Fig.\ref{exact_eg3}.  The simulation  data are generated by selecting  the values of the states at a uniform $5\times 5$ sensor network. At each sensor location, three measurements are taken at time $t=\{0.25,0.75,1\}$, which corresponds to a total of 75 measurements.  For simplicity, the synthetic data $y$ is generated by
\begin{eqnarray*}
y_j=u(\mb{s}_j, t_j; \theta^{\dag})+\xi_j,
\end{eqnarray*}
with $\xi_j\sim N(0, 0.01^2)$. 

 In this example, three hidden layers and 50 neurons per layer are used in DNN.  If the refinement is set to occur, we choose $Q=10$  points to expand the design pool. In order to measure the accuracy of the numerical approximation $\bar{\kappa}$ with respect to the exact solution $\kappa^{\dag}$, we compute the relative error  $rel(k)$   defined as
\begin{eqnarray*}
rel(k)= \frac{\|\bar{\kappa}-\kappa^{\dag}\|_{2}}{\|\kappa^{\dag}\|_{2}},
\end{eqnarray*}
where  $\bar{\kappa}$ is the posterior  mean arising from LDNN or DNN.

\begin{figure}
\begin{center}
         \begin{overpic}[width=0.45\textwidth,trim= 35 0 45 15, clip=true,tics=10]{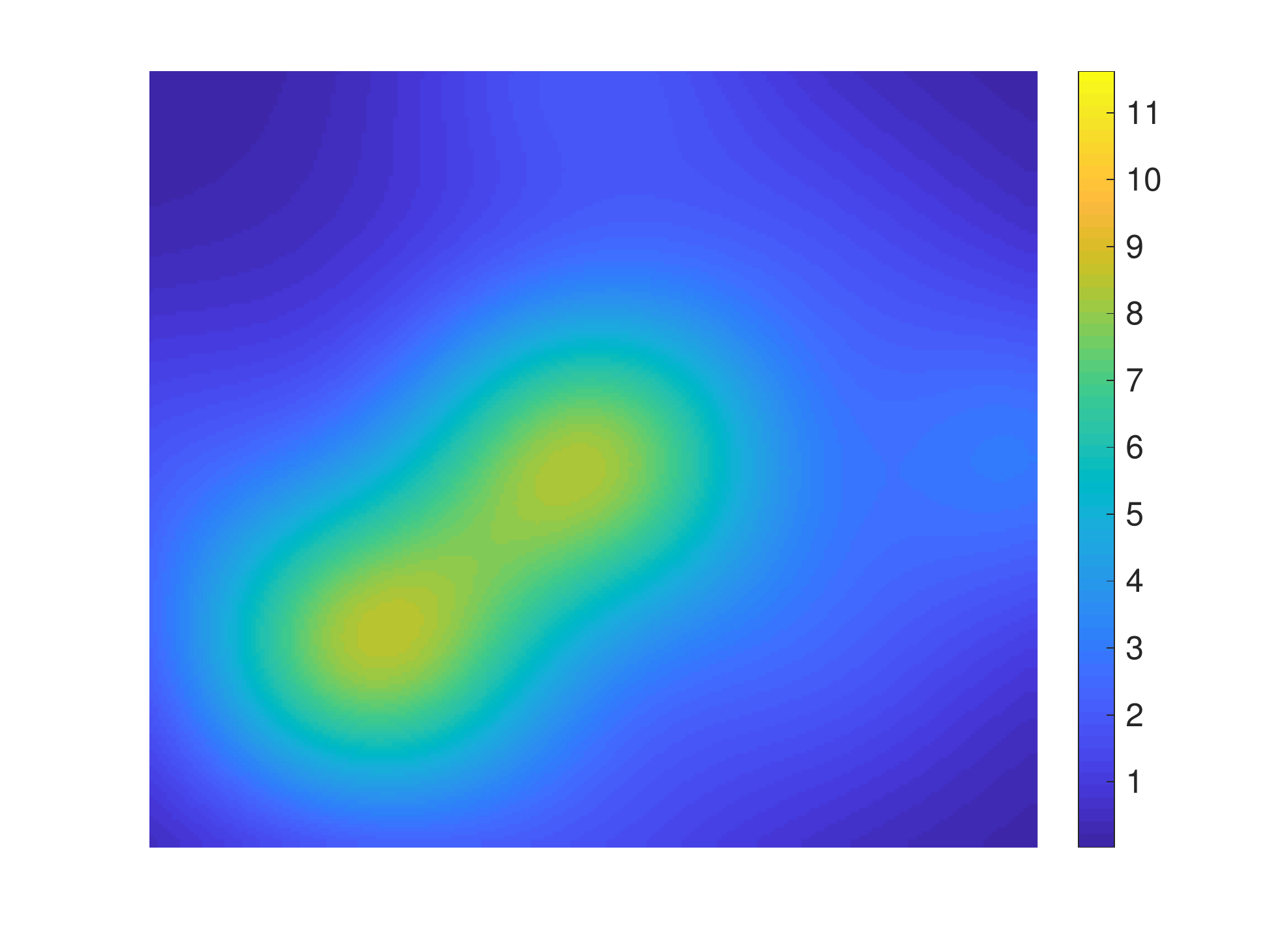}
             \put (35,83) {\footnotesize {\bf DNN, $n_t$=100}}
        \end{overpic}
    \begin{overpic}[width=0.45\textwidth,trim= 35 0 45 15, clip=true,tics=10]{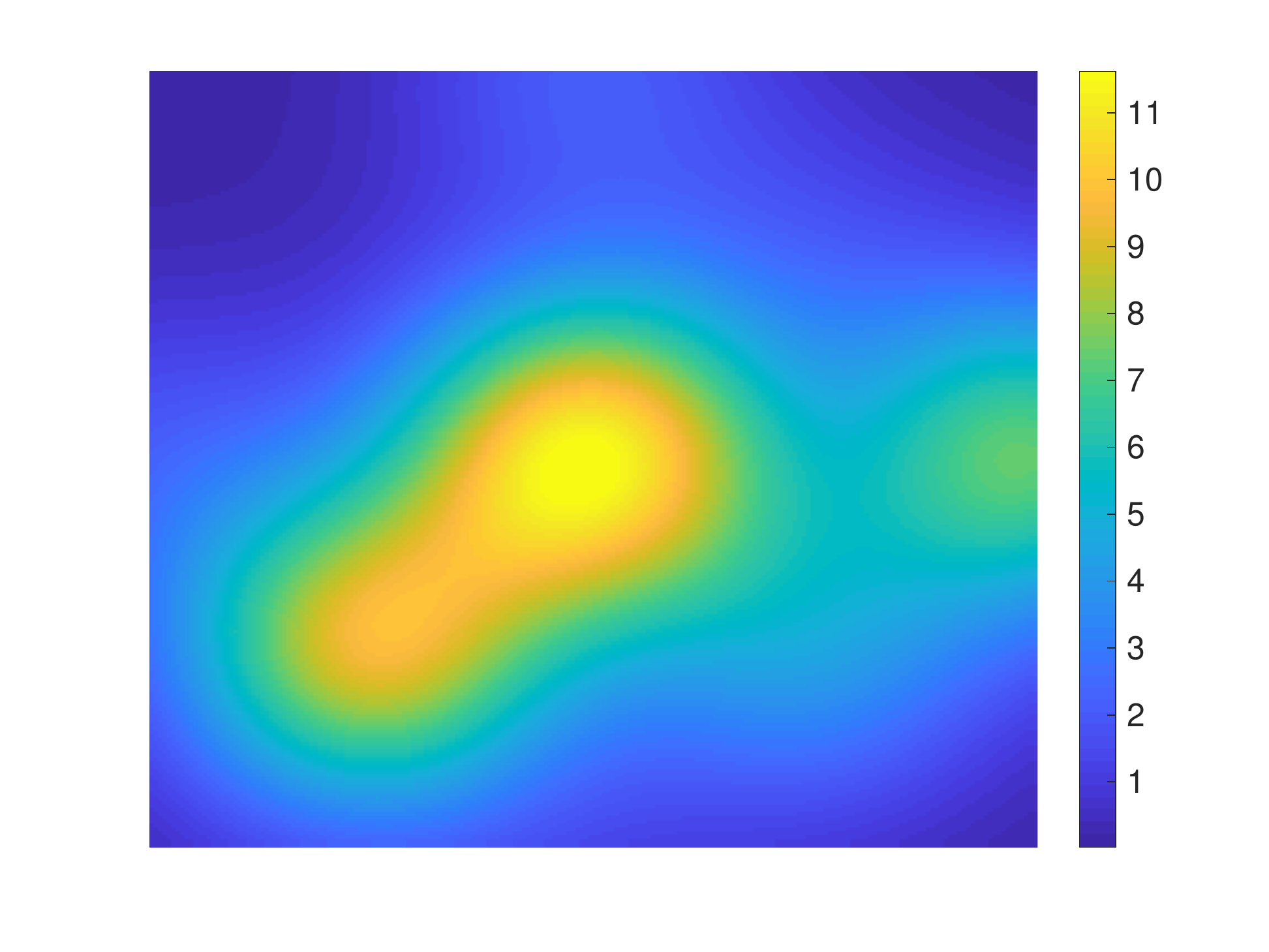}
            \put (35,83) {\footnotesize {\bf LDNN, $n_t$=100}}
      \end{overpic}
      \begin{overpic}[width=0.45\textwidth,trim= 30 0 45 15, clip=true,tics=10]{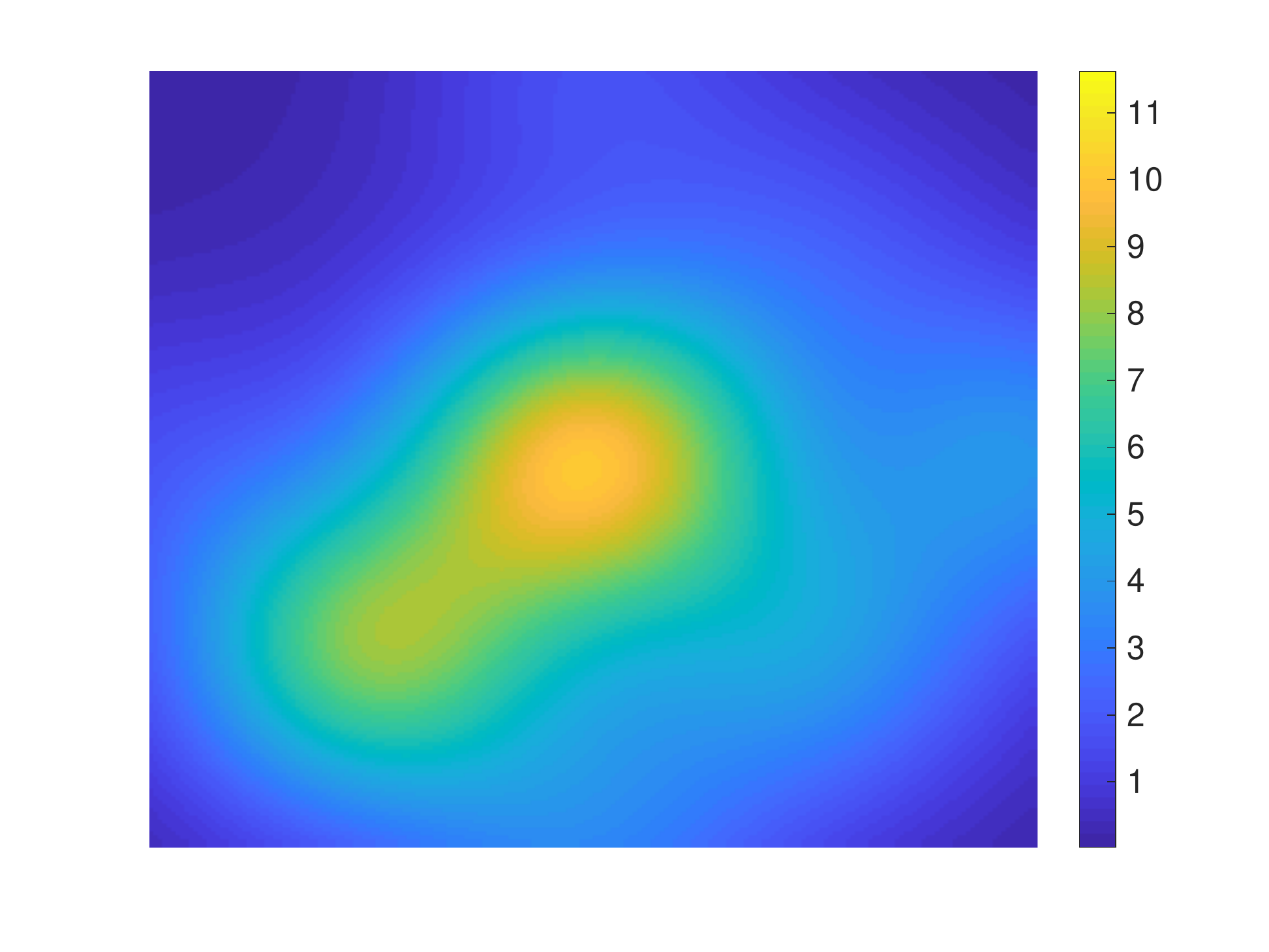}
           \put (35,83) {\footnotesize {\bf DNN,$n_t$=500}}
        \end{overpic}
     \begin{overpic}[width=0.45\textwidth,trim= 30 0 45 15, clip=true,tics=10]{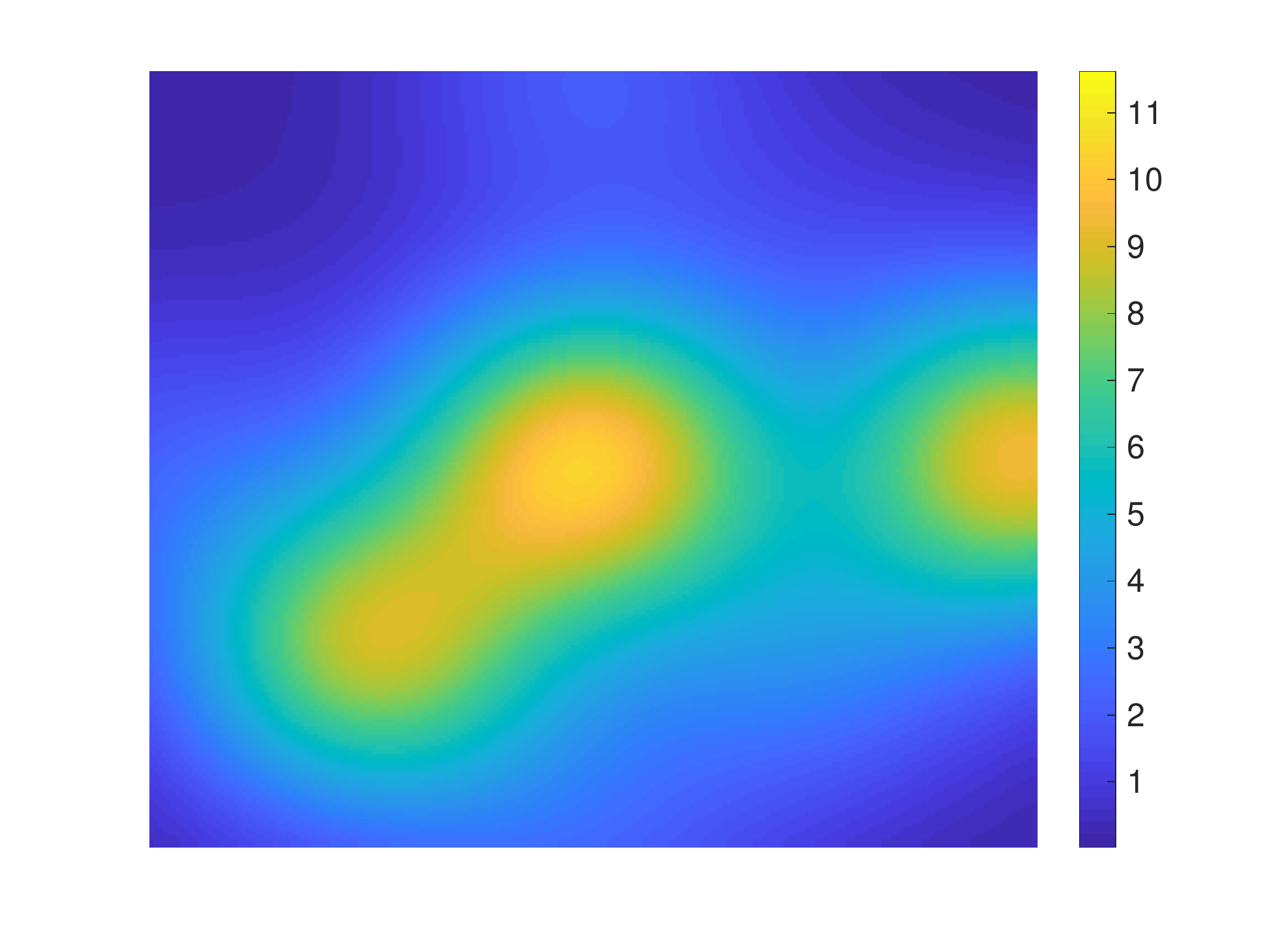}
       \put (35,83) {\footnotesize {\bf LDNN, $n_t$=500}}
        \end{overpic}
     \end{center}
\caption{Posterior mean arising from DNN and LDNN.  From top to bottom, the number of the training set $n_t$ is 100 and 500, respectively. }\label{resLDNN_eg3}
  \end{figure}

  \begin{figure}
\begin{center}
         \begin{overpic}[width=0.5\textwidth,trim= 10 0 35 15, clip=true,tics=10]{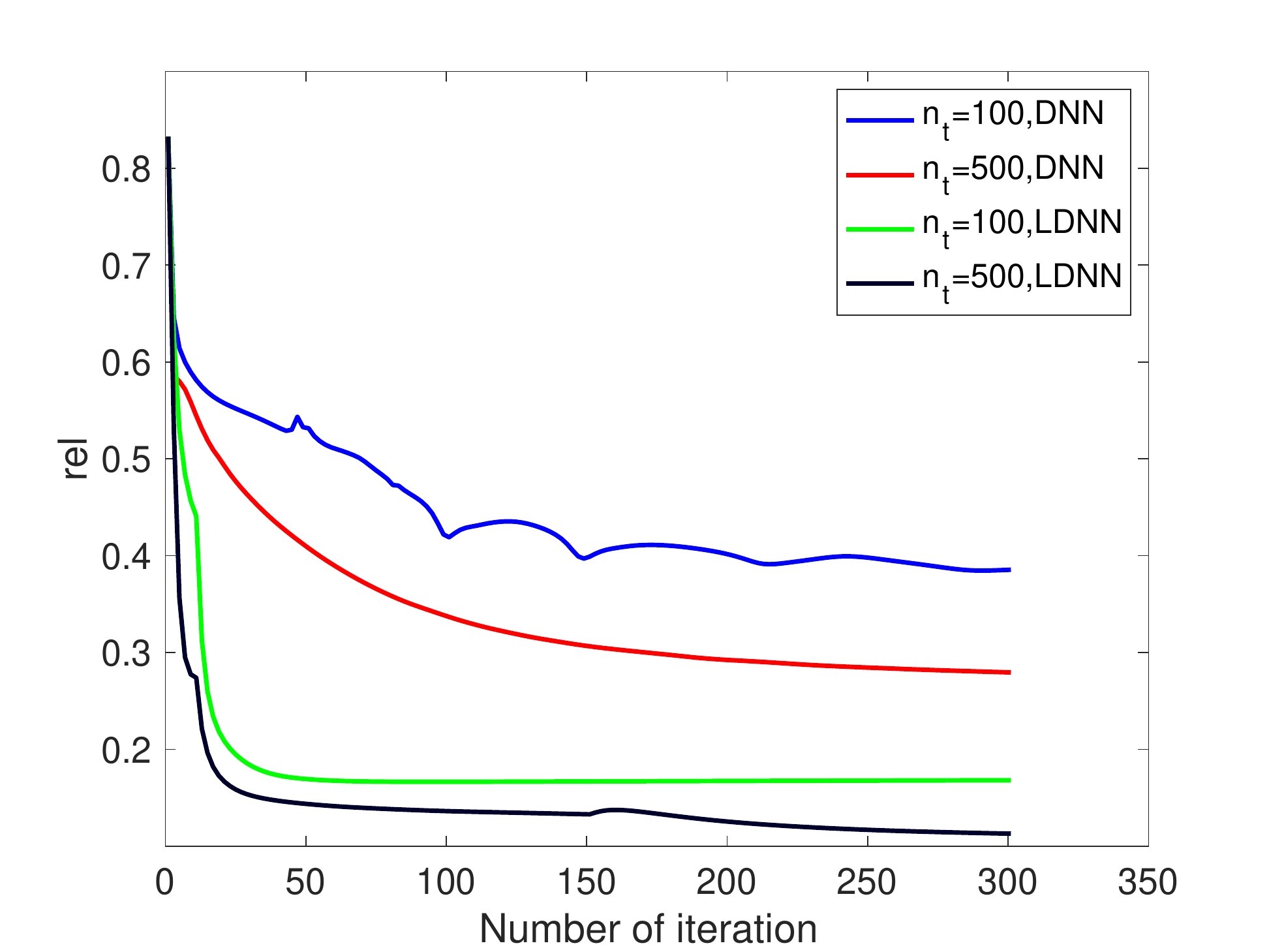}
        \end{overpic}
     \end{center}
\caption{ The accuracy error $rel$ vs. training iteration of different algorithms. }\label{err_eg3}
  \end{figure}
  
   \begin{table}[t]
      \caption{Example 3. Computational times, in seconds, given by three different methods. }\label{eg3_time}
  \centering
  \fontsize{6}{12}\selectfont
  \begin{threeparttable}
    \begin{tabular}{l cccccc}
  \toprule
 & \multicolumn{2}{c}{Offline}&\multicolumn{2}{c}{Online}\cr
\cmidrule(lr){2-3} \cmidrule(lr){4-5}

  \multirow{1}{*}{Method}  &$\text{$\#$ of model eval.}$&CPU(s) &$\text{$\#$ of model eval.}$&CPU(s)     &\multirow{1}{*}{Total time(s)}&\multirow{1}{*}{rel}\cr
  \midrule
   SVGD                           & $-$       & $-$         & $500\times 300 $          & $\sim$3710     &$\sim$ 3710   & $-$  \cr
 
DNN                       & 100         & 37.87           &   $-$                               & 57.97        & 95.84    & 0.3968\cr
DNN                       & 500         & 235.06         &  $-$                               & 57.43       &  292.49   & 0.2941\cr
LDNN                    & 100         & 37.87            &  50                                 &  78.68       & 116.55   & 0.1732\cr
LDNN                    & 500        & 235.06          &  80                                  & 111.26       & 346.32    &0.1155\cr
      \bottomrule
      \end{tabular}
    \end{threeparttable}

\end{table}

In order to verify the accuracy of our proposed algorithm, we compute the posterior
mean using $N=500$ particles arising from the DNN and the LDNN model with different sizes of the training dataset $n_t=\{100, 500\}$.  The numerical results obtained by DNN are shown in the left column of Fig. \ref{resLDNN_eg3}. The corresponding relative errors $rel(k)$ with respect to the training iteration are shown in Fig. \ref{err_eg3}. It can be seen that the numerical results obtained by DNN  results a poor estimate. The corresponding results obtained by LDNN are also shown in Figs. \ref{resLDNN_eg3} and \ref{err_eg3}. It is clearly shown that the LDNN approach results in a very good approximation to the exact solution. Even with a smaller $n_t=100$, the LDNN approach admits a rather accurate result. The total number of high-fidelity model evaluations and the total computational time for LDNN and DNN  are summarized in Table \ref{eg3_time}. Again, the online computational time required by LDNN and DNN is only a small fraction of that by the conventional SVGD. Here, we also use the cost of  $N_{iter}\times N$ high-fidelity model evaluations to represent the CPU time of the original SVGD approach.  It can also be seen from these figures that the LDNN offers a significant improvement in the accuracy, but does not significantly increase the computation time compared to the prior-based DNN approach.
\section{Summary} \label{sec:summary}
The standard SVGD requires the gradient information of the target distribution and cannot be applied when the gradient is unavailable or too expensive to evaluate.  In this work, we introduced a new framework  to address this challenge. The new approaches introduce local approximations of the forward model into the SVGD and refine these approximations incrementally.  One specific scheme, based on deep neural networks (DNN), has been described. The new scheme does not require evaluation of the gradient of the target distribution, thus expanding the application of the SVGD.  The numerical results show that the local approximation can produce accurate results using dramatically fewer evaluations of the forward model.  The proposed methods also provide an emulator that approximates the true posterior density, which can be  employed in further statistical analyses.

While our work focus on vanilla SVGD with local DNN approximation, we believe that other inferential tasks based on  Stein operators can benefit from these developments. Prime candidates include other Stein operators\cite{Liu2017riemannian} and Ensemble Kalman sampler\cite{Garbuno2020interacting}. Integrating subsampling techniques, e.g., \cite{Li2020stochastic}, into the SVGD computation is another promising direction, as the result could more closely mimic standard SVGDs while offering comparable computational savings. Also, the methodology proposed in this paper  may be combined with parameter reduction techniques, as in \cite{Chen2019projected}.   Finally, the proposed framework may be used in conjunction with other reduced-order models, in dynamic data assimilation problems, and for other applications.


\begin{thebibliography}{10}

\bibitem{Bishop2006pattern}
Christopher~M Bishop.
\newblock {\em Pattern recognition and machine learning}.
\newblock springer, 2006.

\bibitem{Kaipio+Somersalo2005}
J.~P. Kaipio and E.~Somersalo.
\newblock {\em Statistical and {C}omputational {I}nverse {P}roblems}, volume
  160.
\newblock Springer, 2005.

\bibitem{Stuart2010}
A.~M. Stuart.
\newblock Inverse problems: a {B}ayesian perspective.
\newblock {\em Acta Numerica}, 19(1):451--559, 2010.

\bibitem{Gelman2013bayesian}
A.~Gelman, J.~B. Carlin, H.~S. Stern, D.~B. Dunson, A.~Vehtari, and D.~B.
  Rubin.
\newblock {\em Bayesian data analysis}.
\newblock CRC press, 2013.

\bibitem{Brooks2011}
S.~Brooks, A.~Gelman, G.~L. Jones, and X.~L. Meng, editors.
\newblock {\em Handbook of {M}arkov chain {M}onte {C}arlo}.
\newblock Chapman \& Hall/CRC Handbooks of Modern Statistical Methods. CRC
  Press, Boca Raton, FL, 2011.

\bibitem{Blei2017variational}
D.~M. Blei, A.~Kucukelbir, and J.~D. McAuliffe.
\newblock Variational inference: A review for statisticians.
\newblock {\em Journal of the American statistical Association},
  112(518):859--877, 2017.

\bibitem{Bardsley2014SISC}
J.~M. Bardsley, A.~Solonen, H.~Haario, and M.~Laine.
\newblock Randomize-then-optimize: A method for sampling from posterior
  distributions in nonlinear inverse problems.
\newblock {\em SIAM Journal on Scientific Computing}, 36(4):A1895--A1910, 2014.

\bibitem{Lan2016emulation}
S.~Lan, T.~Bui-Thanh, M.~Christie, and M.~Girolami.
\newblock Emulation of higher-order tensors in manifold monte carlo methods for
  bayesian inverse problems.
\newblock {\em Journal of Computational Physics}, 308:81--101, 2016.

\bibitem{Martin2012SNMC}
J.~Martin, L.~C. Wilcox, C.~Burstedde, and O.~Ghattas.
\newblock A stochastic newton mcmc method for large-scale statistical inverse
  problems with application to seismic inversion.
\newblock {\em SIAM Journal on Scientific Computing}, 34(3):A1460--A1487, 2012.

\bibitem{Conrad2016JASA}
Patrick~R Conrad, Youssef~M Marzouk, Natesh~S Pillai, and Aaron Smith.
\newblock Accelerating asymptotically exact mcmc for computationally intensive
  models via local approximations.
\newblock {\em Journal of the American Statistical Association},
  111(516):1591--1607, 2016.

\bibitem{Frangos+Marzouk+Willcox2010}
M.~Frangos, Y.~Marzouk, K.~Willcox, and B.~van Bloemen~Waanders.
\newblock Surrogate and reduced-order modeling: a comparison of approaches for
  large-scale statistical inverse problems.
\newblock {\em Biegler, L. and Biros, G. and Ghattas, O. and Heinkenschloss, M.
  and Keyes, D. and Mallick, B. and Marzouk, Y. and Tenorio, L. and van Bloemen
  Waanders, B. and Willcox, K. editors, Computational Methods for Large Scale
  Inverse Problems and Uncertainty Quantification, John Wiley \& Sons, UK},
  pages 123--149, 2010.

\bibitem{stuart+teckentrup2016}
A.~M Stuart and A.~Teckentrup.
\newblock Posterior consistency for {G}aussian process approximations of
  {B}ayesian posterior distributions.
\newblock {\em Mathematics of Computation}, 87(310):721--753, 2018.

\bibitem{yan+guo2015}
L.~Yan and L.~Guo.
\newblock Stochastic collocation algorithms using $l_1$-minimization for
  {B}ayesian solution of inverse problems.
\newblock {\em SIAM Journal on Scientific Computing}, 37(3):A1410--A1435, 2015.

\bibitem{Yan+Zhang2017IP}
L.~Yan and Y.X. Zhang.
\newblock Convergence analysis of surrogate-based methods for bayesian inverse
  problems.
\newblock {\em Inverse Problems}, 33(12):125001, 2017.

\bibitem{Yan+Zhou19JCP}
L.~Yan and T.~Zhou.
\newblock Adaptive multi-fidelity polynomial chaos approach to {B}ayesian
  inference in inverse problems.
\newblock {\em Journal of Computational Physics}, 381:110--128, 2019.

\bibitem{Yan+Zhou2019ADNN}
L.~Yan and T.~Zhou.
\newblock An adaptive surrogate modeling based on deep neural networks for
  large-scale bayesian inverse problems.
\newblock {\em Communications in Computational Physics}, 28(5):2180--2205,
  2020.

\bibitem{Chen2019projected}
P.~Chen, K.~Wu, J.~Chen, T.~O'Leary-Roseberry, and O.~Ghattas.
\newblock Projected stein variational newton: A fast and scalable bayesian
  inference method in high dimensions.
\newblock In {\em Advances in Neural Information Processing Systems}, pages
  15130--15139, 2019.

\bibitem{Detommaso2018stein}
G.~Detommaso, T.~Cui, Y.~Marzouk, A.~Spantini, and R.~Scheichl.
\newblock A stein variational newton method.
\newblock In {\em Advances in Neural Information Processing Systems}, pages
  9169--9179, 2018.

\bibitem{Garbuno2020interacting}
A.~Garbuno-Inigo, F.~Hoffmann, W.~Li, and A.~M. Stuart.
\newblock Interacting langevin diffusions: Gradient structure and ensemble
  kalman sampler.
\newblock {\em SIAM Journal on Applied Dynamical Systems}, 19(1):412--441,
  2020.

\bibitem{Han2018stein}
J.~Han and Q.~Liu.
\newblock Stein variational gradient descent without gradient.
\newblock {\em arXiv preprint arXiv:1806.02775}, 2018.

\bibitem{Li2020stochastic}
L.~Li, Y.~Li, J.G. Liu, Z.~Liu, and J.~Lu.
\newblock A stochastic version of stein variational gradient descent for
  efficient sampling.
\newblock {\em Communications in Applied Mathematics and Computational
  Science}, 15(1):37--63, 2020.

\bibitem{Liu2017riemannian}
C.~Liu and J.~Zhu.
\newblock Riemannian stein variational gradient descent for bayesian inference.
\newblock {\em arXiv preprint arXiv:1711.11216}, 2017.

\bibitem{Liu2017stein}
Q.~Liu.
\newblock Stein variational gradient descent as gradient flow.
\newblock In {\em Advances in neural information processing systems}, pages
  3115--3123, 2017.

\bibitem{Liu2016stein}
Q.~Liu and D.~Wang.
\newblock Stein variational gradient descent: A general purpose bayesian
  inference algorithm.
\newblock In {\em Advances in neural information processing systems}, pages
  2378--2386, 2016.

\bibitem{Lu2019scaling}
J.~Lu, Y.~Lu, and J.~Nolen.
\newblock Scaling limit of the stein variational gradient descent: The mean
  field regime.
\newblock {\em SIAM Journal on Mathematical Analysis}, 51(2):648--671, 2019.

\bibitem{Wang2019stein}
D.~Wang, Z.~Tang, C.~Bajaj, and Q.~Liu.
\newblock Stein variational gradient descent with matrix-valued kernels.
\newblock In {\em Advances in neural information processing systems}, pages
  7836--7846, 2019.

\bibitem{Robbins1951stochastic}
H.~Robbins and S.~Monro.
\newblock A stochastic approximation method.
\newblock {\em The annals of mathematical statistics}, pages 400--407, 1951.

\bibitem{Zeiler2012Adadelta}
M.~Zeiler.
\newblock Adadelta: an adaptive learning rate method.
\newblock {\em arXiv preprint arXiv:1212.5701}, 2012.

\bibitem{Wang2020particle}
Y.~Wang, J.~Chen, L.~Kang, and C.~Liu.
\newblock Particle-based energetic variational inference.
\newblock {\em arXiv preprint arXiv:2004.06443}, 2020.

\bibitem{Goodfellow2016DL}
I.~Goodfellow, Y.~Bengio, and A.~Courville.
\newblock {\em Deep learning}.
\newblock MIT press, 2016.

\bibitem{Ramachandran2017}
P.~Ramachandran, B.~Zoph, and Q.~Le.
\newblock Searching for activation functions.
\newblock {\em arXiv preprint arXiv:1710.05941}, 2017.

\bibitem{Tripathy+Bilionis2018JCP}
R.~K. Tripathy and I.~Bilionis.
\newblock Deep {UQ}: Learning deep neural network surrogate models for high
  dimensional uncertainty quantification.
\newblock {\em Journal of Computational Physics}, 375:565--588, 2018.

\bibitem{Bottou2010}
L.~Bottou.
\newblock Large-scale machine learning with stochastic gradient descent.
\newblock In {\em Proceedings of COMPSTAT'2010}, pages 177--186. Springer,
  2010.

\bibitem{Tieleman+Hinton2012lecture}
T.~Tieleman and G.~Hinton.
\newblock Lecture 6.5-rmsprop: Divide the gradient by a running average of its
  recent magnitude.
\newblock {\em COURSERA: Neural networks for machine learning}, 4(2):26--31,
  2012.

\bibitem{Kingma2014Adam}
D.~Kingma and J.~Ba.
\newblock Adam: A method for stochastic optimization.
\newblock {\em arXiv preprint arXiv:1412.6980}, 2014.

\bibitem{Han+Jentzen+E2018PNAS}
J.~Han, A.~Jentzen, and W.~E.
\newblock Solving high-dimensional partial differential equations using deep
  learning.
\newblock {\em Proceedings of the National Academy of Sciences},
  115(34):8505--8510, 2018.

\bibitem{Raissi2019JCP}
M.~Raissi, P.~Perdikaris, and G.~E. Karniadakis.
\newblock Physics-informed neural networks: A deep learning framework for
  solving forward and inverse problems involving nonlinear partial differential
  equations.
\newblock {\em Journal of Computational Physics}, 378:686--707, 2019.

\bibitem{Schwab+Zech2019AA}
C.~Schwab and J.~Zech.
\newblock Deep learning in high dimension: Neural network expression rates for
  generalized polynomial chaos expansions in uq.
\newblock {\em Analysis and Applications}, 17(01):19--55, 2019.

\bibitem{Zhu+Zabaras2018bayesian}
Y.~Zhu and N.~Zabaras.
\newblock Bayesian deep convolutional encoder--decoder networks for surrogate
  modeling and uncertainty quantification.
\newblock {\em Journal of Computational Physics}, 366:415--447, 2018.

\bibitem{Yan+Zhou2019RTO}
L.~Yan and T.~Zhou.
\newblock An acceleration strategy for randomize-then-optimize sampling via
  deep neural networks.
\newblock {\em submitted}, 2020.

\bibitem{Yosinski2014transferable}
J.~Yosinski, J.~Clune, Y.~Bengio, and H.~Lipson.
\newblock How transferable are features in deep neural networks?
\newblock In {\em Advances in neural information processing systems}, pages
  3320--3328, 2014.

\bibitem{Gretton2012kernel}
A.~Gretton, K.~M Borgwardt, M.~J Rasch, B.~Sch{\"o}lkopf, and A.~Smola.
\newblock A kernel two-sample test.
\newblock {\em The Journal of Machine Learning Research}, 13(1):723--773, 2012.

\bibitem{Lin+Xu2007}
Y.~Lin and C.~Xu.
\newblock Finite difference/spectral approximations for the time-fractional
  diffusion equation.
\newblock {\em Journal of Computational Physics}, 225(2):1533--1552, 2007.

\bibitem{Haario+Laine+Mira2006}
H.~Haario, M.~Laine, A.~Mira, and E.~Saksman.
\newblock {DRAM}: efficient adaptive {MCMC}.
\newblock {\em Statistics and Computing}, 16(4):339--354, 2006.

\bibitem{Yan+Zhou2019PCEKI}
L.~Yan and T.~Zhou.
\newblock An adaptive multifidelity pc-based ensemble kalman inversion for
  inverse problems.
\newblock {\em International Journal for Uncertainty Quantification},
  9(3):205--220, 2019.

\end{thebibliography}
\end{document}